\documentclass[10pt]{article}
\newcommand{\hidetag}[1]{}
\usepackage{graphics}
\usepackage{graphicx}
\usepackage{color}
\usepackage{epsfig}
\usepackage{amsfonts}

\usepackage{amsmath}
\usepackage{amstext}
\usepackage{amsopn}
\usepackage{amsbsy}
\usepackage{amscd}
\usepackage{amsxtra}
\usepackage{amsthm}
\usepackage{graphicx}
\usepackage{epsfig}
\usepackage{pstricks,multido,pst-node}

\usepackage{fancyhdr}
\numberwithin{equation}{section}

\newcommand{\phie}{\phi_\epsilon}
\newcommand{\ds}{\displaystyle}
\newtheorem{theorem}{Theorem}[section]
\newtheorem{lemma}[theorem]{Lemma}
\newtheorem{proposition}[theorem]{Proposition}
\newtheorem{corollary}[theorem]{Corollary}
\newtheorem{remark}[theorem]{Remark}

\newcommand{\BE}{\begin{equation}}
\newcommand{\BEN}{\begin{equation*}}
\newcommand{\EE}{\end{equation}}
\newcommand{\EEN}{\end{equation*}}
\newcommand{\BL}{\begin{lemma}}
\newcommand{\EL}{\end{lemma}}
\newcommand{\BT}{\begin{theorem}}
\newcommand{\ET}{\end{theorem}}
\newcommand{\BP}{\begin{proposition}}
\newcommand{\EP}{\end{proposition}}
\newcommand{\BC}{\begin{corollary}}
\newcommand{\EC}{\end{corollary}}
\newcommand{\BR}{\begin{remark}}
\newcommand{\ER}{\end{remark}}

\title{Multiple solutions of steady-state Poisson-Nernst-Planck equations with steric effects }

\author{Tai-Chia Lin \thanks{Institute of
Applied Mathematical Sciences, Center for Advanced Study in Theoretical Sciences (CASTS), National Taiwan University, No.1, Sec.4, Roosevelt Road, Taipei 106,
Taiwan, email: {\tt tclin@math.ntu.edu.tw}} \and Bob Eisenberg\thanks{Department of Molecular Biophysics $\&$ Physiology Rush Medical Center, 1653 West Congress, Parkway, Chicago, IL 60612, USA, email: {\tt beisenbe@rush.edu} } }

\date{}

\begin{document}

%\allowdisplaybreaks
\maketitle
%-----------------------------------------------------------------------------------------------------------
%% Abstract
%-----------------------------------------------------------------------------------------------------------
\begin{abstract}
Experiments measuring currents through single protein channels show unstable currents. Channels switch between 'open' or 'closed' states in a spontaneous stochastic process called gating. Currents are either (nearly) zero or at a definite level, characteristic of each type of protein, independent of time, once the channel is open. The steady state Poisson-Nernst-Planck equations with steric effects (PNP-steric equations) describe steady current through the open channel quite well, in a wide variety of conditions. Here we study the existence of multiple solutions of steady state PNP-steric equations to see if they themselves, without modification or augmentation, can describe two levels of current. We prove that there are two steady state solutions of PNP-steric equations for (a) three types of ion species (two types of cations and one type of anion) with a positive constant permanent charge, and (b) four types of ion species (two types of cations and their counter-ions) with a constant permanent charge but no sign condition. The excess currents (due to steric effects) associated with these two steady state solutions are derived and expressed as two distinct formulas. Our results indicate that PNP-steric equations may become a useful model to study spontaneous gating of ion channels. Spontaneous gating is thought to involve small structural changes in the channel protein that perhaps produce large changes in the profiles of free energy that determine ion flow. Gating is known to be modulated by external structures. Both can be included in future extensions of our present analysis.
\end{abstract}


%\begin{AMS}76A05, 76M99, 65C30\end{AMS}
%-----------------------------------------------------------------------------------------------------------

\section{Introduction}
\medskip
\noindent

The Poisson-Nernst-Planck (PNP) equations, a well-known model of ion transport, play a crucial role in the study of many physical and biological phenomena (cf.~\cite{BCEJ97, BTA04, BMPS02, B_ip11, CK_itn05, EC93, EL07, E08, J_tpm10, LZ12, MRS90, S84, WZCX_sr12}). Such an important model can be represented by
	
\BE\label{pnp1}
\left\{ \begin{array}{ll}
   & \frac{\partial {{c}_{i}}}{\partial t}+\nabla \cdot J_{i}^{PNP}=0,\quad i=1,\cdots, N\,,  \\
   & \\
   &  -J_{i}^{PNP}={{D}_{i}}\left( \nabla {{c}_{i}}+\frac{{{z}_{i}}e}{{{k}_{B}}T}{{c}_{i}}\nabla\phi\right)\,, \\
   & \\
   &-\nabla \cdot (\varepsilon \nabla \phi )={{\rho }_{0}}+\sum\limits_{i=1}^{N}{{{z}_{i}}e{{c}_{i}}}  \\
\end{array} \right.
\EE
where $N$ is the number of ion species, $c_i$ is the distribution function, $J_{i}^{PNP}$ is the flux density, $D_i$ is the diffusion constant, and $z_i$ is the valence of the $i$th ion species, respectively. Besides, $\phi$ is the electrostatic potential, $\varepsilon$ is the dielectric constant, $\rho_0$ is the permanent (fixed) charge density of the system, $k_B$ is the Boltzmann constant, $T$ is the absolute temperature and $e$ is the elementary charge. Due to ionic sizes, steric repulsion may appear in crowded ions of several biological systems like DNAs, ribosomes and ion channels. When ions are crowded in a narrow channel, the PNP equations become unreliable because the ion-size effect becomes important, but the PNP equations represent ions as point particles without size (cf.~\cite{AGAN02, BSK11, Er11, F-jpc10, HP-jpc10, KSD-PRL10, Lb2, VVB-jcp10}).

To include ion size effects, Eisenberg and Liu modified PNP equations into a complicated system of differential-integral equations with singular integrals that simulate successfully the selectivity of important types of calcium and sodium ion channels (cf.~\cite{HEL2}). However, the singular integrals form an extremely singular kernel because of the Lennard-Jones (LJ) potential. Numerical efficiency and theoretical analysis disappear when forced to deal with such singularities (cf.~\cite{Er1, HEL}). To simplify the model, we truncate the (spatial) frequency range of the LJ potential, find a simpler energy functional from the leading order terms of the energy expansion with suitable scales. We derive the Poisson-Nernst-Planck equations with steric effects called PNP-steric equations (cf.~\cite{LE-cms13})
\begin{align}
\label{eqn1}
&\frac{\partial {{c}_{i}}}{\partial t}+\nabla \cdot {{J}_{i}}=0, i=1,\cdots ,N\,,\\
\label{eqn2}
&-\nabla \cdot (\varepsilon \nabla \phi )={{\rho}_{0}}+\sum\limits_{i=1}^{N}{{{z}_{i}}e{{c}_{i}}}\,,
\end{align}
where flux ${{J}_{i}}$ is
\BE\label{eqn3}	
{{J}_{i}}=-{{D}_{i}}\nabla {{c}_{i}}-\frac{{{D}_{i}}{{c}_{i}}}{{{k}_{B}}T}{{z}_{i}}e\nabla \phi -\frac{{{D}_{i}}{{c}_{i}}}{{{k}_{B}}T}\sum\limits_{j=1}^{N}{{{g}_{ij}}\nabla {{c}_{j}}}\,,
\EE
and ${{g}_{ij}}={{g}_{ji}}\sim\epsilon_{ij} {{\left( {{a}_{i}}+{{a}_{j}}\right)}^{12}}$ is a nonnegative constant depending on ion radii ${{a}_{i}}$, ${{a}_{j}}$ and the energy coupling constant $\epsilon_{ij}$ of the $i$-th and $j$-th species ions, respectively (cf.~\cite{HLLE}). Note that equations (\ref{eqn1})-(\ref{eqn3}) can be regarded as a system of reaction-diffusion equations with nonlinear cross-diffusion terms being similar to~\cite{BFPS-sjma10}. Amazingly, these equations are an effective model to simulate the selectivity of ion channels (cf.~\cite{HLLE}).

Comparing (\ref{eqn3}) with $J_{i}^{PNP}$ in (\ref{pnp1}), the excess flux $J_{i}^{ex}={{J}_{i}}-J_{i}^{PNP}$ due to steric effects of ion species $i$ is
$$
-J_{i}^{ex}=\frac{1}{{{k}_{B}}T}{{D}_{i}}{{c}_{i}}\nabla\mu_{i}^{ex}\quad\hbox{ and }\quad\mu_{i}^{ex}=\sum\limits_{j=1}^{N}{{{g}_{ij}}{{c}_{j}}}
$$
where $\mu_{i}^{ex}=\sum\limits_{j=1}^{N}{{{g}_{ij}}{{c}_{j}}}$ is the excess chemical potential of ion species $i$ due to steric effects. Consequently, the excess current ${{I}^{ex}}=\sum\limits_{i=1}^{N}{{{z}_{i}}eJ_{i}^{ex}}$
due to steric effects becomes
\BE\label{I-ex1}
{{I}^{ex}}=-\sum\limits_{i,j=1}^{N}\,\frac{z_i\,e}{{{k}_{B}}T}{{D}_{i}}\,{{g}_{ij}}\,{{c}_{i}}\nabla{{c}_{j}}\,.
\EE
We shall use the formula (\ref{I-ex1}) to calculate the excess currents for multiple solutions of the 1D steady-state PNP-steric equations. We are motivated by the hope--but cannot dare expect--that one solution will correspond to a closed state and the other to an open state, as found in experiments~\cite{E-ws96} and in simulations~\cite{KLTME-pb2013}. Of course, the current measured through the open state corresponds to the total current, not just the excess currents.

The existence of multiple steady (equilibrium) states is important to study transitions between such states which may be related to the gating (switching between open and closed states) and selectivity of ion channels. Multiple steady states can be investigated by finding multiple solutions of the 1D steady-state PNP equations for two types of ion species with three regions of piecewise constant permanent charge under the assumption that the Debye number is large~\cite{EL07}. More general theorems related to multiple solutions of the 1D steady-state PNP equations involving multiple types of ions with multiple regions of piecewise constant permanent charge are discussed in~\cite{ws-Liu-jde09}. With only a constant permanent charge, there is only a unique solution of the 1D steady-state PNP equations for multiple types of ions~\cite{MRS90, WHWH_pre2014}.
Instead of the 1D steady-state PNP equations, here we study multiple solutions of the 1D steady-state PNP-steric equations with spatially constant permanent charges.

For simplicity, we consider domain as a 1D interval $(-1,1)$ for (\ref{eqn1})-(\ref {eqn3}) and set  ${{J}_{i}}=0$,  $i=1,\cdots ,N$ to get the steady-state PNP-steric equations. Then by (\ref {eqn3}),
\[\frac{d}{dx}\left( \ln {{c}_{i}}+\frac{{{z}_{i}}e}{{{k}_{B}}T}\phi +\frac{1}{{{k}_{B}}T}\sum\limits_{j=1}^{N}{{{g}_{ij}}{{c}_{j}}} \right)=0 \quad\hbox{for}\: x\in (-1,1)\,, i=1,\cdots,N\,,\]
which can be satisfied if
\BE\label{eqn3-1}
\ln {{c}_{i}}+\frac{{{z}_{i}}e}{{{k}_{B}}T}\phi +\frac{1}{{{k}_{B}}T}\sum\limits_{j=1}^{N}{{{g}_{ij}}{{c}_{j}}}=0\quad\hbox{ for }\:i=1,\cdots,N\,,
\EE
holds true. Let $\tilde{\phi }=\frac{e}{{{k}_{B}}T}\phi $ and ${{\tilde{g}}_{ij}}=\frac{1}{{{k}_{B}}T}{{g}_{ij}}$ for  $i,j=1,\cdots ,N$. Then (\ref{eqn2}) and (\ref{eqn3-1}) can be transformed into
\BE\label{eqn3-2}	
\ln{{c}_{i}}+{{z}_{i}}\tilde{\phi }+\sum\limits_{j=1}^{N}{{{{\tilde{g}}}_{ij}}{{c}_{j}}}=0\quad \text{ for }\;i=1,\cdots ,N,
\EE
and
\BE\label{eqn3-3}
-\tilde{\varepsilon }{{\tilde{\phi }}_{xx}}={{\tilde{\rho }}_{0}}+\sum\limits_{i=1}^{N}{{{z}_{i}}{{c}_{i}}}\quad\hbox{ for }\: x\in (-1,1)\,,
\EE
where  $\tilde{\varepsilon }=\frac{{{k}_{B}}T}{{{e}^{2}}}\varepsilon $ and  ${{\tilde{\rho }}_{0}}=\frac{1}{e}{{\rho }_{0}}$. For notational convenience, we may remove tilde ($\sim$) and denote (\ref{eqn3-2}) and (\ref{eqn3-3}) as
\BE\label{eqn3-4}	
\ln {{c}_{i}}+{{z}_{i}}\phi +\sum\limits_{j=1}^{N}{{{g}_{ij}}{{c}_{j}}}=0\quad \text{ for }\;i=1,\cdots ,N,
\EE
and
\BE\label{eqn-1.9}
-\varepsilon{{\phi }_{xx}}={{\rho }_{0}}+\sum\limits_{i=1}^{N}{{{z}_{i}}{{c}_{i}}}\quad\hbox{ for }\: x\in (-1,1)\,.
\EE
Equations like (\ref{eqn3-4}) have been used to interpret bioelectric phenomena in many papers since they were adopted by Hodgkin, Huxley, and Cole (cf.~\cite{C_jp68, H_tn02}). Here we consider the following boundary condition given by
\BE\label{bdc-Rbn-0}
\phi(1)+\eta_\varepsilon\phi'(1)=\phi_0(1)
\quad\hbox{and}\quad
\phi(-1)-\eta_\varepsilon\phi'(-1)=\phi_0(-1)\,,\EE
where $\phi_0(1), \phi_0(-1)$ are constants and $\eta_\varepsilon$ is a
non-negative constant. Here ${{\phi }_{0}}\left( \pm 1 \right)$ and $\phi \left( \pm 1 \right)$ are the extrachannel and intrachannel electrostatic potentials at the channel boundaries, respectively. The coefficient ${{\eta }_{\varepsilon }}\sim \frac{{{\varepsilon }_{0}}}{{{\varepsilon }_{m}}}$ is governed by the ratio of  ${{\varepsilon }_{0}}$ the dielectric constant of the electrolyte solution and ${{\varepsilon }_{m}}$ the dielectric constant of the membrane (cf~\cite{ZBL}). Note that (\ref{bdc-Rbn-0}) is of the Robin boundary condition if $\eta_\varepsilon>0$; and of the Dirichlet boundary condition if $\eta_\varepsilon=0$. The Robin boundary condition includes polarization (e.g. dielectric) charges in the bath and/or electrodes which the Dirichlet boundary condition does not. Such charges, induced by and dependent on the electric field play a prominent role in the art of real experiments, because they are important determinants of the background noise and stability of high speed recordings. The theoretical reasons for these practical realities have not been investigated to the best of our knowledge.

As $N=2$, the existence, uniqueness and the solution's asymptotic behavior of (\ref{eqn3-4})-(\ref{bdc-Rbn-0}) are investigated under non-symmetry breaking condition $0\le {{g}_{12}}=g_{21}\le \sqrt{{{g}_{11}}{{g}_{22}}}$ which implies that solution $\left( {{c}_{1}},{{c}_{2}} \right)$ of (\ref{eqn3-4}) is uniquely determined by $\phi$ (cf.~\cite{LLL-pt13}). Hence (\ref{eqn3-4}) and (\ref{eqn-1.9}) can be reduced to a single differential equation of $\phi$. However, as the symmetry breaking condition ${{g}_{12}}=g_{21}> \sqrt{{{g}_{11}}{{g}_{22}}}$ holds true, solution $\left( {{c}_{1}},{{c}_{2}} \right)$ of (\ref{eqn3-4}) may not be uniquely determined by $\phi$. In Section~\ref{ss-solu-2s}, we introduce new variables $\xi, \Sigma $ and transform (\ref{eqn3-4}) into a quadratic polynomial which can be solved precisely to get explicit formulas and represent two branches of solution curves. Using these explicit formulas, we can then define biological conductance (for that condition) as the biologists do and perform the comparison using formulas like (\ref{iA-3-int})-(\ref{iB-3BN}). Note that the symbol $g$ is used for conductance (units siemens) in biology and this is not equivalent to our $g_{ij}$. In this paper, we want to study multiple solutions of (\ref{eqn3-4})-(\ref{bdc-Rbn-0}) for the cases of $N=3, 4$, and $g_{12}=g_{21}, g_{34}=g_{43}$ sufficiently large such that symmetry breaking condition ${{g}_{12}}=g_{21}> \sqrt{{{g}_{11}}{{g}_{22}}}, {{g}_{34}}=g_{43}> \sqrt{{{g}_{33}}{{g}_{44}}}$ holds true.

\subsection{Main Results}
\ \ \
System (\ref{eqn3-4}) can be regarded as a coupled system of algebraic equations. Because ${{g}_{ij}}=0$ for $i,j=1,\cdots,N$, a solution of system (\ref{eqn3-4}) can be expressed as ${{c}_{i}}={{e}^{-{{z}_{i}}\phi }}$ for $i=1,\cdots,N$. However, it seems impossible to solve system (\ref{eqn3-4}) explicitly for the general case of ${{g}_{ij}}>0$ for $i,j=1,\cdots,N$. To overcome such difficulty, we may set $N=2$, $z_2=-z_1=q\geq 1$, $g_{11}=g_{22}=g>0$, and introduce new variables $\xi ={{c}_{1}}{{c}_{2}}$ and $\Sigma ={{c}_{1}}+{{c}_{2}}$. Then (\ref{eqn3-4}) can be transformed into a quadratic polynomial that can be solved explicitly (see Section~\ref{ss-solu-2s}). For $g_{12}=g_{21}=z$ large (see Theorem~\ref{pro-cn-cp} in Section~\ref{ss-solu-2s}), system (\ref{eqn3-4}) has two branches of solutions $\left({{c}_{1}},{{c}_{2}}\right)=\left({{c}_{1}}\left({{\Sigma }_{{{A}_{1}}}}\left(\phi\right)\right),{{c}_{2}}\left({{\Sigma}_{{{A}_{1}}}}\left(\phi\right)\right)\right)$ and $\left({{c}_{1}},{{c}_{2}} \right)=\left({{c}_{1}}\left({{\Sigma }_{{{B}_{1}}}}\left(\phi\right) \right),{{c}_{2}}\left({{\Sigma }_{{{B}_{1}}}}\left(\phi\right)\right)\right)$ such that $(c_1-c_2)\circ\Sigma_{A_1}:[-\phi_{A,c},\infty)\to\mathbb{R}$ and $(c_1-c_2)\circ\Sigma_{B_1}:(-\infty,\phi_{A,c}]\to\mathbb{R}$ are monotone increasing functions to $\phi$, where $\phi_{A,c}>0$ is a constant, $\Sigma_{A_1}$ and $\Sigma_{B_1}$ are two functions satisfying
\begin{align}\notag & (c_1-c_2)\circ\Sigma_{A_1}(-\phi_{A,c})=(c_1-c_2)(\Sigma_c)>0\,,\\
&
(c_1-c_2)\circ\Sigma_{B_1}(\phi_{A,c})=-(c_1-c_2)(\Sigma_c)<0\,,\notag\end{align}
$$\lim_{\phi\to\infty}(c_1-c_2)\circ\Sigma_{A_1}(\phi)=\infty\quad\hbox{and}\: \lim_{\phi\to-\infty}(c_1-c_2)\circ\Sigma_{B_1}(\phi)=-\infty\,.$$
Here $\circ$ denotes the function $(c_1 -c_2)$ acting on the function $\Sigma_{A_1}(\phi)$, i.e., the function composition and $g_c$ is the positive constant defined in Proposition~\ref{rho_c}. Besides, $\phi_{A,c}$ satisfies $\phi_{A,c}\to +\infty$ and $(c_1-c_2)(\Sigma_c)\to 0$ as $z\to +\infty$ and $g>0$ is fixed. Hence (\ref{eqn3-4}) and (\ref{eqn-1.9}) can be decomposed into two differential equations like (\ref{eqn3-5-A}) and (\ref{eqn3-5-B}) but they can not have uniformly bounded solutions to $\varepsilon>0$ (see Lemma~\ref{ub-solu}). This fact motivates us to add one extra species $c_3$ and assume that $N=3$, $g_{12}=g_{21}=z$ is sufficiently large, $g_{11}=g_{22}=g>0$, $z_2=-z_1=q\geq 1, z_3>0$, $g_{i3}=g_{3i}=0$, $i=1,2,3$ (which implies $c_3=e^{-z_3\phi}$). Then (\ref{eqn3-4}) and (\ref{eqn-1.9}) may be reduced to two differential equations (\ref{eqn3-5-A}) and (\ref{eqn3-5-B}) having uniformly bounded solutions, respectively. This may provide multiple solutions of (\ref{eqn3-4})-(\ref{bdc-Rbn-0}).

Natural biological solutions always contain at least three species (sodium, potassium, and chloride, and usually calcium). Experiments are often done, however, with just two species (say sodium chloride) along with traces of hydrogen ion, and perhaps other contaminants. Gating occurs in simplified unnatural situations and so we hope to study  mathematical solutions in corresponding situations in a separate paper.

Now we state the main result of this paper as follows:
\BT\label{thm1.1}
Let $N=3$, $z_2=-z_1=q\geq 1, z_3>0$ and $\rho_0>0$ be a constant. Assume that $g_{11}=g_{22}=g>0$ is fixed and $g_{i3}=g_{3i}=0$ for $i=1,2,3$. Then as $g_{12}=g_{21}=z>0$ is sufficiently large, the system of equations (\ref{eqn3-4})-(\ref{bdc-Rbn-0}) has two uniformly bounded (to $\varepsilon$) solutions $\phi_{\varepsilon}^{A}$ and $\phi_{\varepsilon }^{B}$ such that $\phi_{\varepsilon}^{A}\left( x \right)\to {{\phi}_{{{A}_{1}},0}}$ and $\phi_{\varepsilon }^{B}\left( x \right)\to {{\phi }_{{{B}_{1}},0}}$ for $x\in (-1,1)$ as $\varepsilon\rightarrow 0$, where ${{\phi}_{{{A}_{1}},0}}$ and ${{\phi }_{{{B}_{1}},0}}$ are two distinct constants.
\ET
\noindent In most of the "cation" (e.g., sodium, potassium, and calcium) channels, $\rho_0$ is a negative number. There are regions ('rings') of negative charge and some channels (sodium channel DEKA) have a ring of positive charge as well. Here we assume the positive sign of $\rho_0$ which may produce the values ${{\phi}_{{{A}_{1}},0}}$ and ${{\phi }_{{{B}_{1}},0}}$ (see Figure~\ref{fig4} in Section~\ref{mul-solu-3s-1}), and the proof of Theorem~\ref{thm1.1} is given in Section~\ref{mul-solu-3s-1}.

To remove the sign condition on $\rho_0$, we may consider four ion species composed of two cations and counterions (like the mixture of ${\rm N{{a}^{+}},{{Ca}^{+2}},C{{l}^{-}}}$ and ${\rm CO_3^{-2}}$) and study multiple solutions of (\ref{eqn3-4})-(\ref{bdc-Rbn-0}) with $N=4$, $z_2=-z_1=q_1\geq 1$, $z_4=-z_3=q_2\geq 1$, $g_{11}=g_{22}=g>0$, and $g_{33}=g_{44}=\tilde{g}>0$. Using the assumption ${{g}_{ij}}={{g}_{ji}}=0$ for $i=1,2$ and $j=3,4$, we may decompose system (\ref{eqn3-4}) with $N=4$ into two independent systems having the same form as (\ref{eqn3-4}) with $N=2$. Hence Theorem~\ref{pro-cn-cp} (in Section~\ref{ss-solu-2s}) implies that as $g_{12}=g_{21}=z$ and $g_{34}=g_{43}=\tilde{z}>0$ sufficiently large, system (\ref{eqn3-4}) has four branches of solutions
\begin{eqnarray*}
&&\left({{c}_{1}},{{c}_{2}}\right)=\left({{c}_{1}}\left({{\Sigma }_{{{A}_{1}}}}\left(\phi\right)\right),{{c}_{2}}\left({{\Sigma}_{{{A}_{1}}}}\left(\phi\right)\right)\right)\,,\quad  \left({{c}_{1}},{{c}_{2}} \right)=\left({{c}_{1}}\left({{\Sigma }_{{{B}_{1}}}}\left(\phi\right) \right),{{c}_{2}}\left({{\Sigma }_{{{B}_{1}}}}\left(\phi\right)\right)\right)\,, \\
&&\left({{c}_{3}},{{c}_{4}}\right)=\left({{c}_{3}}\left({{\Sigma }_{{{M}_{1}}}}\left(\phi\right)\right),{{c}_{4}}\left({{\Sigma}_{{{M}_{1}}}}\left(\phi\right)\right)\right)\,,\quad  \left({{c}_{3}},{{c}_{4}} \right)=\left({{c}_{3}}\left({{\Sigma }_{{{N}_{1}}}}\left(\phi\right) \right),{{c}_{4}}\left({{\Sigma}_{{{N}_{1}}}}\left(\phi\right)\right)\right)\,,
\end{eqnarray*}
such that
$(c_1-c_2)\circ\Sigma_{A_1}:[-\phi_{A,c},\infty)\to\mathbb{R}$, $(c_1-c_2)\circ\Sigma_{B_1}:(-\infty,\phi_{A,c}]\to\mathbb{R}$,
$(c_3-c_4)\circ\Sigma_{N_1}:[-\phi_{M,c},\infty)\to\mathbb{R}$ and $(c_3-c_4)\circ\Sigma_{M_1}:(-\infty,\phi_{M,c}]\to\mathbb{R}$, are monotone increasing functions of $\phi$, where $\phi_{A,c}, \phi_{M,c}>0$ are constants, $\Sigma_{A_1}$, $\Sigma_{B_1}$, $\Sigma_{M_1}$ and $\Sigma_{N_1}$ are functions satisfying
\begin{align*} & (c_1-c_2)\circ\Sigma_{A_1}(-\phi_{A,c})\,, (c_3-c_4) \circ\Sigma_{N_1}(-\phi_{M,c})>0\,,\\
& (c_1-c_2)\circ\Sigma_{B_1}(\phi_{A,c})\,, (c_3-c_4) \circ\Sigma_{M_1}(\phi_{M,c})<0\,, \\
& \lim_{\phi\to\infty}(c_1-c_2)\circ\Sigma_{A_1}(\phi)= \lim_{\phi\to\infty}(c_3-c_4)\circ\Sigma_{N_1}(\phi)=\infty\,, \\
& \lim_{\phi\to-\infty}(c_1-c_2)\circ\Sigma_{B_1}(\phi)= \lim_{\phi\to-\infty}(c_3-c_4)\circ\Sigma_{M_1}(\phi)=-\infty\,.
\end{align*}
Here $\circ$ denotes function composition. Moreover, $\phi_{A,c}, \phi_{M,c}\to +\infty$ and $(c_1-c_2) \circ\Sigma_{A_1}(-\phi_{A,c}), (c_1-c_2)\circ\Sigma_{B_1}(\phi_{A,c}), (c_3-c_4) \circ\Sigma_{M_1}(\phi_{M,c})$ and $(c_3-c_4) \circ\Sigma_{N_1}(-\phi_{M,c})$ tend to zero as $z, \tilde{z}\to +\infty$ and $g, \tilde{g}>0$ are fixed.

Without loss of generality, we may assume $\phi_{M,c}<\phi_{A,c}$. Then the graphs of functions $(c_1-c_2)\circ\Sigma_{A_1}$ and $(c_4-c_3) \circ\Sigma_{M_1}$ may intersect at $\phi=\phi_{A_1,0}$ as $z$ and $\tilde{z}$ sufficiently large (see Figure~5 in Section~\ref{mul-solu-3s-2}). Similarly, the graphs of functions $(c_1-c_2)\circ\Sigma_{B_1}$ and $(c_4-c_3) \circ\Sigma_{N_1}$ may intersect at $\phi=\phi_{B_1,0}$ as $z$ and $\tilde{z}$ sufficiently large. Hence (\ref{eqn3-4}) and (\ref{eqn-1.9}) may be reduced to two differential equations with the same forms as (\ref{eqn3-5-A}) and (\ref{eqn3-5-B}) having uniformly bounded solutions, respectively. This may provide the following result for multiple solutions of (\ref{eqn3-4})-(\ref{bdc-Rbn-0}).
\BT\label{thm1.2}
Let $N=4$, $z_2=-z_1=q_1\geq 1$, $z_4=-z_3=q_2\geq 1$ and $\rho_0\neq 0$ be a constant. Assume that $g_{11}=g_{22}=g>0$, $g_{33}=g_{44}=\tilde{g}>0$ are fixed and ${{g}_{ij}}={{g}_{ji}}=0$ for $i=1,2$ and $j=3,4$. Then as $g_{12}=g_{21}=z>0$ and $g_{34}=g_{43}=\tilde{z}>0$ are sufficiently large, the system of equations (\ref{eqn3-4})-(\ref{bdc-Rbn-0}) has two uniformly bounded (to $\varepsilon$) solutions $\phi_{\varepsilon}^{A}$ and $\phi_{\varepsilon }^{B}$ such that $\phi_{\varepsilon}^{A}\left( x \right)\to {{\phi}_{{{A}_{1}},0}}$ and $\phi_{\varepsilon }^{B}\left( x \right)\to {{\phi }_{{{B}_{1}},0}}$ for $x\in (-1,1)$ as $\varepsilon\rightarrow 0$, where ${{\phi}_{{{A}_{1}},0}}$ and ${{\phi }_{{{B}_{1}},0}}$ are two distinct constants.
\ET
\noindent  The proof of Theorem~\ref{thm1.2} is given in Section~\ref{mul-solu-3s-2}.

For solutions $\phi_{\varepsilon }^{A}$ and $\phi_{\varepsilon }^{B}$, the corresponding excess currents defined in (\ref{I-ex1}) may be denoted as $I_{A}^{ex}$ and $I_{B}^{ex}$, respectively. Under the same hypotheses of Theorem~\ref{thm1.1} for three ion species, we may use the new variable $\Sigma $ to derive the following formulas (see Section~\ref{ex-cu1}):
\BE\label{iA-3-int}
\begin{array}{rll}
&\ds\int_{{{x}_{1}}}^{{{x}_{2}}}{I_{A}^{ex}dx} \\
& \\
&=q\,e\,\ds\int_{{{\Sigma}^A_{1}}}^{{{\Sigma}^A_{2}}}\,\frac{{{D}_{2}}-{{D}_{1}}}{2}\left\{(1-q)-q\left[ g\Sigma +\left( {{g}^{2}}-{{z}^{2}} \right){{e}^{-(g+z)\Sigma }} \right]\right\}d\Sigma \\
& \\
& \hspace{0.3cm} -q\,e\,\ds\int_{{{\Sigma}^A_{1}}}^{{{\Sigma}^A_{2}}}\,\frac{{{D}_{1}}+{{D}_{2}}}{2\sqrt{{{\Sigma }^{2}}-4{{e}^{-(g+z)\Sigma }}}}\left\{(1-q)\Sigma-q\,g{{\Sigma }^{2}}+(g+z)\left[2-q\,\left( {{g}}-{{z}} \right)\Sigma\right] {{e}^{-(g+z)\Sigma }} \right\}d\Sigma\,,
\end{array}
\EE
and
\BE\label{iB-3-int}
\begin{array}{rll}
&\ds\int_{{{x}_{1}}}^{{{x}_{2}}}{I_{B}^{ex}dx} \\ & \\ &=q\,e\,\ds\int_{{{\Sigma}^B_{1}}}^{{{\Sigma}^B_{2}}}{}
\frac{{{D}_{2}}-{{D}_{1}}}{2}\left\{(1-q)-q\left[ g\Sigma +\left( {{g}^{2}}-{{z}^{2}} \right){{e}^{-(g+z)\Sigma }} \right]\right\}d\Sigma \\
& \\ & \hspace{0.3cm} +q\,e\ds\int_{{{\Sigma}^B_{1}}}^{{{\Sigma }^B_{2}}}{}\frac{{{D}_{1}}+{{D}_{2}}}{2\sqrt{{{\Sigma }^{2}}-4{{e}^{-(g+z)\Sigma }}}}\left\{(1-q)\Sigma-q\,g{{\Sigma }^{2}}+(g+z)\left[2-q\,\left( {{g}}-{{z}} \right)\Sigma\right] {{e}^{-(g+z)\Sigma }} \right\}d\Sigma\,,
\end{array}
\EE
for $-1<x_1<x_2<1$, where ${{\Sigma}^A_{j}}={{\Sigma }_{{{A}_{1}}}}\left( {{\phi}^A_{\varepsilon}(x_j)} \right)$ and ${{\Sigma}^B_{j}}={{\Sigma }_{{{B}_{1}}}}\left( {{\phi}^B_{\varepsilon}(x_j)} \right)$ for $j=1,2$. From (\ref{iA-3-int}) and (\ref{iB-3-int}), it is clear that the difference between $I_{A}^{ex}$ and $I_{B}^{ex}$ which may give various ion flows related to currents observed in channels as they switch (i.e., gate) from one level of current to another.

The method of Section~\ref{ex-cu1} can be generalized to four ion species with the same hypotheses of Theorem~\ref{thm1.2}. As for (\ref{iA-3-int}) and (\ref{iB-3-int}), we may derive (see Section~\ref{ex-cu2})
\BE\label{iA-3AM}
\begin{array}{rll}
&\ds\int_{{{x}_{1}}}^{{{x}_{2}}}{I_{A, M}^{ex}dx} = \ds\int_{{{x}_{1}}}^{{{x}_{2}}}{I_{A}^{ex}+I_{M}^{ex}dx} \\
&=q_1\,e\,\ds\int_{{{\Sigma}^A_{1}}}^{{{\Sigma}^A_{2}}}\,\frac{{{D}_{2}}-{{D}_{1}}}{2}\left\{(1-q_1)-q_1\left[ g\Sigma +\left( {{g}^{2}}-{{z}^{2}} \right){{e}^{-(g+z)\Sigma }} \right]\right\}d\Sigma \\
& \\
& \hspace{0.3cm} -q_1\,e\,\ds\int_{{{\Sigma}^A_{1}}}^{{{\Sigma}^A_{2}}}\,\frac{{{D}_{1}}+{{D}_{2}}}{2\sqrt{{{\Sigma }^{2}}-4{{e}^{-(g+z)\Sigma }}}}\left\{(1-q_1)\Sigma-q_1\,g{{\Sigma }^{2}}+(g+z)\left[2-q_1\,\left( {{g}}-{{z}} \right)\Sigma\right] {{e}^{-(g+z)\Sigma }} \right\}d\Sigma \\
& \\
& \hspace{0.3cm} +q_2\,e\,\ds\int_{{{\Sigma}^M_{1}}}^{{{\Sigma}^M_{2}}}\,\frac{{{D}_{4}}-{{D}_{3}}}{2}\left\{(1-q_2)-q_2\left[ \tilde{g}\Sigma +\left( {{\tilde{g}}^{2}}-{{\tilde{z}}^{2}} \right){{e}^{-(\tilde{g}+\tilde{z})\Sigma }} \right]\right\}d\Sigma \\
& \\
& \hspace{0.3cm} -q_2\,e\,\ds\int_{{{\Sigma}^M_{1}}}^{{{\Sigma}^M_{2}}}\,\frac{{{D}_{3}}+{{D}_{4}}}{2\sqrt{{{\Sigma }^{2}}-4{{e}^{-(\tilde{g}+\tilde{z})\Sigma }}}}\left\{(1-q_2)\Sigma-q_2\,\tilde{g}{{\Sigma }^{2}}+(\tilde{g}+\tilde{z})\left[2-q_2\,\left( {{\tilde{g}}}-{{\tilde{z}}} \right)\Sigma\right] {{e}^{-(\tilde{g}+\tilde{z})\Sigma }} \right\}d\Sigma\,,
\end{array}
\EE
and
\BE\label{iB-3BN}
\begin{array}{rll}
&\ds\int_{{{x}_{1}}}^{{{x}_{2}}}{I_{B,N}^{ex}dx}=\ds\int_{{{x}_{1}}}^{{{x}_{2}}}{I_{B}^{ex}+I_{N}^{ex}dx} \\ &=q_1\,e\,\ds\int_{{{\Sigma}^B_{1}}}^{{{\Sigma}^B_{2}}}{}
\frac{{{D}_{2}}-{{D}_{1}}}{2}\left\{(1-q_1)-q_1\left[ g\Sigma +\left( {{g}^{2}}-{{z}^{2}} \right){{e}^{-(g+z)\Sigma }} \right]\right\}d\Sigma \\
& \\ & \hspace{0.3cm} +q_1\,e\,\ds\int_{{{\Sigma}^B_{1}}}^{{{\Sigma }^B_{2}}}{}\frac{{{D}_{1}}+{{D}_{2}}}{2\sqrt{{{\Sigma }^{2}}-4{{e}^{-(g+z)\Sigma }}}}\left\{(1-q_1)\Sigma-q_1\,g{{\Sigma }^{2}}+(g+z)\left[2-q_1\,\left( {{g}}-{{z}} \right)\Sigma\right] {{e}^{-(g+z)\Sigma }} \right\}d\Sigma \\
& \\ & \hspace{0.3cm} +q_2\,e\,\ds\int_{{{\Sigma}^N_{1}}}^{{{\Sigma}^N_{2}}}\,\frac{{{D}_{4}}-{{D}_{3}}}{2}\left\{(1-q_2)-q_2\left[ \tilde{g}\Sigma +\left( {{\tilde{g}}^{2}}-{{\tilde{z}}^{2}} \right){{e}^{-(\tilde{g}+\tilde{z})\Sigma }} \right]\right\}d\Sigma \\
& \\
& \hspace{0.3cm} +q_2\,e\,\ds\int_{{{\Sigma}^N_{1}}}^{{{\Sigma}^N_{2}}}\,\frac{{{D}_{3}}+{{D}_{4}}}{2\sqrt{{{\Sigma }^{2}}-4{{e}^{-(\tilde{g}+\tilde{z})\Sigma }}}}\left\{(1-q_2)\Sigma-q_2\,\tilde{g}{{\Sigma }^{2}}+(\tilde{g}+\tilde{z})\left[2-q_2\,\left( {{\tilde{g}}}-{{\tilde{z}}} \right)\Sigma\right] {{e}^{-(\tilde{g}+\tilde{z})\Sigma }} \right\}d\Sigma\,,
\end{array}
\EE
where ${{\Sigma}^A_{j}}={{\Sigma }_{{{A}_{1}}}}\left(\phi_{\varepsilon }^{A}\left( {{x}_{j}} \right) \right)$, ${{\Sigma}^M_{j}}={{\Sigma }_{{{M}_{1}}}}\left(\phi_{\varepsilon }^{A}\left( {{x}_{j}} \right) \right)$, ${{\Sigma}^B_{j}}={{\Sigma }_{{{B}_{1}}}}\left(\phi_{\varepsilon }^{B}\left( {{x}_{j}} \right) \right)$, and ${{\Sigma}^N_{j}}={{\Sigma }_{{{N}_{1}}}}\left(\phi_{\varepsilon }^{B}\left( {{x}_{j}} \right) \right)$  for $j=1,2$.
The difference between $I_{A,M}^{ex}$ and $I_{B,N}^{ex}$ may also give various ion flows related to currents observed in channels as they switch (i.e., gate) from one level of current to another.

The rest of this paper is organized as follows: We may solve system (\ref{eqn3-4}) of algebraic equations with $N=2$, $z_2=-z_1=q\geq 1$ and $g_{11}=g_{22}>0$ in Section~\ref{ss-solu-2s}. Theorem~\ref{thm1.1} and~\ref{thm1.2} are proven in Section~\ref{mul-solu-3s}. The proofs of Lemma~\ref{bd-solu} and~\ref{ub-solu} are given in Section~\ref{sec4}, and formulas (\ref{iA-3-int})-(\ref{iB-3BN}) are derived in Section~\ref{ex-cu}.

\section{Solutions of (\ref{eqn3-4}) with $N=2$, $z_2=-z_1=q\geq 1$ and $g_{11}=g_{22}$}\label{ss-solu-2s}
\ \ \ In this section, we study equation (\ref{eqn3-4}) with $N=2$, $z_2=-z_1=q\geq 1$ and $g_{11}=g_{22}=g$ which can be  denoted as follows:
\begin{align} \label{eqn1-5-4} \left(\ln c_1-q\,\phi\right) + \left(g\,c_1 + z\,c_2\right) &=0\,, \\
\label{eqn2-5-4} \left(\ln c_2+q\,\phi\right) + \left(g\,c_2 + z\,c_1\right) &=0\,,
\end{align}
where $z=g_{12}$ and $g=g_{11}=g_{22}$ are positive constants. Physically, $g_{ij}\sim\epsilon_{ij} \left(a_i+a_j\right)^{12}$, where $a_i$ is the ion radius of $i$-th ion species with concentration $c_i$, and $\epsilon_{ij}>0$ is the energy coupling constant between $i$-th and $j$-th ion species for $i=1,2$. Note that (\ref{eqn1-5-4}) and (\ref{eqn2-5-4}) are formulated as a system of algebraic equations. We want to solve these equations and get solutions for $\left(c_1,c_2\right)$ as a function of $\phi$. Adding (\ref{eqn1-5-4}) and (\ref{eqn2-5-4}), we get
\BE\label{2.1}
\ln \left( {{c}_{1}}{{c}_{2}} \right)+\left( g+z\right)\left( {{c}_{1}}+{{c}_{2}} \right)=0\,.
\EE

Now we introduce new variables as follows:
\[\xi ={{c}_{1}}{{c}_{2}}\quad \text{ and }\quad \Sigma={{c}_{1}}+{{c}_{2}}.\]
Multiplying  $\sigma$ by ${{c}_{1}}$ , we get a quadratic polynomial of $c_1$ as follows:
$$
\Sigma {{c}_{1}}=c_{1}^{2}+\xi
$$
which gives  ${{c}_{1}}=\frac{\Sigma \pm \sqrt{{{\Sigma }^{2}}-4\xi }}{2}$ and hence by  ${{c}_{1}}{{c}_{2}}=\xi$,  $\left( {{c}_{1}},{{c}_{2}} \right)$ can be expressed as
\BE\label{eqn15420}
\begin{array}{ll}
   & \left({{c}_{1}},{{c}_{2}}\right)=\left(\frac{\Sigma +\sqrt{{{\Sigma }^{2}}-4\xi }}{2},\frac{\Sigma -\sqrt{{{\Sigma }^{2}}-4\xi }}{2} \right)\,,
  \\
   &\text{or}  \\
&\left({{c}_{1}},{{c}_{2}}\right)=\left( \frac{\Sigma -\sqrt{{{\Sigma }^{2}}-4\xi }}{2},\frac{\Sigma +\sqrt{{{\Sigma }^{2}}-4\xi }}{2} \right)\,,  \\
\end{array}
\EE
for $\Sigma\geq 2\sqrt{\xi}>0$. Moreover, (\ref{2.1}) can be transformed into  $\ln \xi =-\left( g+z \right)\Sigma $
i.e.
\BE\label{2.2}
\xi ={{e}^{-\left( g+z \right)\Sigma }}.
\EE
Hence the solution $(c_1,c_2)$ of (\ref{eqn1-5-4}) and (\ref{eqn2-5-4}) may be described by two curves $A$ and $B$ parameterized by the total concentration $\Sigma$ and denoted as
\BE\label{cvA}
A=\left\{ ({{c}_{1}},{{c}_{2}})=\left( \frac{\Sigma +\sqrt{{{\Sigma }^{2}}-4{{e}^{-\left( g+z \right)\Sigma }}}}{2},\frac{\Sigma -\sqrt{{{\Sigma }^{2}}-4{{e}^{-\left( g+z \right)\Sigma }}}}{2} \right):\Sigma \ge {{\Sigma }_{z}}>0 \right\}\,,
\EE
and
\BE\label{cvB}
B=\left\{ ({{c}_{1}},{{c}_{2}})=\left( \frac{\Sigma -\sqrt{{{\Sigma }^{2}}-4{{e}^{-\left( g+z \right)\Sigma }}}}{2},\frac{\Sigma +\sqrt{{{\Sigma }^{2}}-4{{e}^{-\left( g+z \right)\Sigma }}}}{2} \right):\Sigma \ge {{\Sigma }_{z}}>0 \right\}\,.
\EE
Here $\Sigma_z>0$ a critical total concentration is the unique positive solution of ${{\Sigma }^{2}}=4{{e}^{-\left( g+z \right)\Sigma }}$ such that concentrations $c_1$ and $c_2$ are equal to $\frac{1}{2}\Sigma_z$ as the total concentration $\Sigma=\Sigma_z$. Then
\BE\label{cn-cp}
{{c}_{1}}-{{c}_{2}}=\left\{ \begin{array}{rrrll}
   {} & \sqrt{{{\Sigma }^{2}}-4{{e}^{-\left( g+z \right)\Sigma }}} & \text{ on }\;A,  \\
   {} & -\sqrt{{{\Sigma }^{2}}-4{{e}^{-\left( g+z \right)\Sigma }}} & \text{on}\;B.  \\
\end{array} \right.
\EE

Take (\ref{cvA}) and (\ref{cvB}) into (\ref{eqn1-5-4}), and let $\phi_A=\phi$ on curve $A$, and $\phi_B=\phi$ on curve $B$, respectively. Then
\BE\label{ph1}
\left\{\begin{array}{rlll} & q\,{{\phi }_{A}}(\Sigma )=\ln \left[ \frac{1}{2}\left( \Sigma +\sqrt{{{\Sigma }^{2}}-4{{e}^{-\left( g+z \right)\Sigma }}} \right) \right]+\frac{g+z}{2}\Sigma +\frac{g-z}{2}\sqrt{{{\Sigma }^{2}}-4{{e}^{-\left( g+z \right)\Sigma }}}\,,\\
& \\
& q\,{{\phi }_{B}}(\Sigma )=\ln \left[ \frac{1}{2}\left( \Sigma -\sqrt{{{\Sigma }^{2}}-4{{e}^{-\left( g+z \right)\Sigma }}} \right) \right]+\frac{g+z}{2}\Sigma +\frac{z-g}{2}\sqrt{{{\Sigma }^{2}}-4{{e}^{-\left( g+z \right)\Sigma }}}\,,
\end{array}\right. \EE
for $\Sigma\geq\Sigma_z$. Consequently,
\BE\label{sym}
\phi_A+\phi_B=0\,,\EE
\BE\label{dphia}
q\,\frac{d{{\phi }_{A}}}{d\Sigma }=\frac{(1+g\Sigma ){{e}^{(g+z)\Sigma }}+{{g}^{2}}-{{z}^{2}}}{{{e}^{(g+z)\Sigma }}\sqrt{{{\Sigma }^{2}}-4{{e}^{-(g+z)\Sigma }}}}\,,
\EE
and
\BE\label{dphib}
q\,\frac{d{{\phi }_{B}}}{d\Sigma }=-\frac{(1+g\Sigma ){{e}^{(g+z)\Sigma }}+{{g}^{2}}-{{z}^{2}}}{{{e}^{(g+z)\Sigma }}\sqrt{{{\Sigma }^{2}}-4{{e}^{-(g+z)\Sigma }}}}\,,
\EE
when the total concentration $\Sigma$ is larger than $\Sigma_z$. Note that curve $A$ and $B$ are joined only at a single point  $({{c}_{1}},{{c}_{2}})=\left( \frac{1}{2}{{\Sigma }_{z}},\frac{1}{2}{{\Sigma }_{z}} \right)$ which is located only at $\Sigma=\Sigma_z$. Moreover, ${{\phi }_{A}}\left( {{\Sigma }_{z}} \right)={{\phi }_{B}}\left( {{\Sigma }_{z}} \right)=0$ and $\left( {{c}_{1}}-{{c}_{2}} \right)\left( {{\Sigma }_{z}} \right)=0$.

Suppose $0<z\le g$. Then $\phi_A$ and $\phi_B$ can be regarded as one variable $\phi$ and $c_1-c_2$ may become a strictly monotone increasing function of $\phi$. The result is stated as follows:
\BP\label{mon}
Suppose $0<z\le g$. Then $\Sigma=\Sigma(\phi)$ can be a single-valued function of $\phi$ with domain being the entire space $\mathbb{R}$ and range  $\left[ {{\Sigma }_{z}},\infty  \right)$ such that $\Sigma (0)={{\Sigma }_{z}}$,
\BE\label{inv-ro1}
\left\{ \begin{array}{lll}
   &{{\phi }_{A}}(\Sigma(\phi))=\phi &\text{ if }\phi \ge 0\,,  \\
&  & \\
   &{{\phi }_{B}}(\Sigma(\phi))=\phi &\text{ if }\phi \le 0\,,  \\
\end{array} \right.
\EE
and ${{c}_{1}}-{{c}_{2}}=\left( {{c}_{1}}-{{c}_{2}} \right)\left( \Sigma \left( \phi  \right) \right)$ is a strictly monotone increasing function of $\phi$ from $-\infty$ to $\infty$.
\EP
\begin{proof}
Suppose $0<z\le g$. Then by (\ref{dphia}) and (\ref{dphib}), we have
\BE\label{mon2}
\frac{d}{d\Sigma }{{\phi }_{A}}\left( \Sigma \right)>0\quad\hbox{ and}\quad \frac{d}{d\Sigma }{{\phi }_{B}}\left( \Sigma \right)<0\quad\hbox{for}\:\Sigma\geq\Sigma_z\,.
\EE
Here we have used $0<z\le g$. Thus $\phi_A(\Sigma)>0$ and $\phi_B(\Sigma)<0$ for $\Sigma>\Sigma_z$. Besides, the range of $\phi_A$ is $[0,\infty)$ and the range of $\phi_B$ is $(-\infty,0]$. Note that ${{\phi }_{A}}({{\Sigma }_{z}})={{\phi }_{B}}({{\Sigma }_{z}})=0$. We may combine $\phi_A$ and $\phi_B$ as one variable $\phi$ (see Figure~\ref{fig1}) defined as follows:
\[\left\{ \begin{array}{lll}
   &\phi={{\phi }_{A}}(\Sigma)\ge 0 &\text{ on } A\,,  \\
&  & \\
   &\phi={{\phi }_{B}}(\Sigma)\le 0 &\text{ on } B\,.  \\
\end{array} \right.\]

%\vspace*{4cm}
%\begin{figure}[htp]
%	\centering
%	\hspace*{1cm}
%	\includegraphics[scale=0.45]{rho_phi_1.eps}
%	\vspace*{-5cm}
%   \caption{$\theta = (\rho_z,0)$ in $(\rho,\phi)$ coordinates}
%    \label{fig1}
%\end{figure}
%\newpage
\vspace*{-2cm}
\begin{figure}[htp]
	\centering
	\includegraphics[scale=0.35, bb=0 0 800 570]{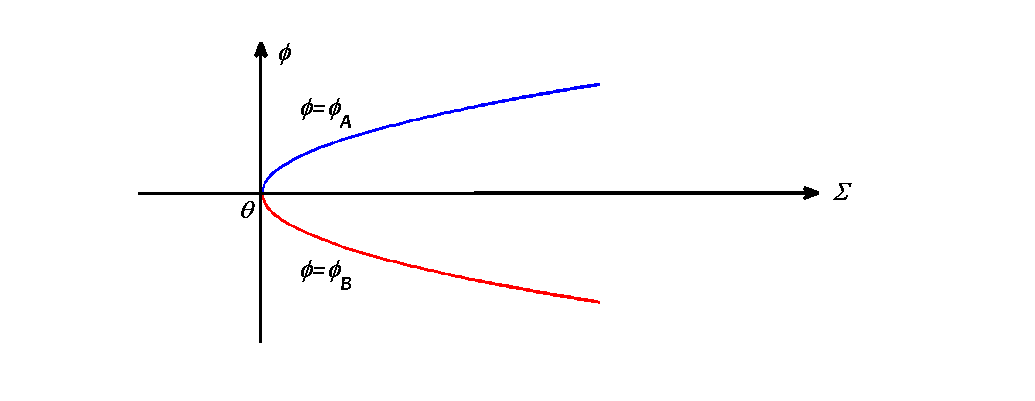}
	\vspace*{-0.5cm}
    \caption{$\theta = (\Sigma_z,0)$ in $(\Sigma,\phi)$ coordinates}
    \label{fig1}
\end{figure}
%\vspace*{-0.5cm}
\noindent Hence by (\ref{mon2}) and inverse function theorem, $\Sigma$ can be denoted as $\Sigma=\Sigma(\phi)$ and become a single-valued function of $\phi$ with domain being the entire space $\mathbb{R}$ and range  $\left[{{\Sigma }_{z}},\infty\right)$ such that $\Sigma(0)=\Sigma_z$ and (\ref{inv-ro1}) hold true. The derivative of $\Sigma$ with respect to $\phi$ is
\BE\label{drho-phi}
\frac{d\Sigma }{d\phi }=\frac{1}{\frac{d\phi }{d\Sigma }}=\left\{ \begin{array}{rllll}
   & q\,\frac{{{e}^{(g+z)\Sigma }}\sqrt{{{\Sigma }^{2}}-4{{e}^{-(g+z)\Sigma }}}}{(1+g\Sigma ){{e}^{(g+z)\Sigma }}+{{g}^{2}}-{{z}^{2}}} &\text{ if } & \phi \ge 0\,,\\
& & & \\
& -q\,\frac{{{e}^{(g+z)\Sigma }}\sqrt{{{\Sigma }^{2}}-4{{e}^{-(g+z)\Sigma }}}}{(1+g\Sigma ){{e}^{(g+z)\Sigma }}+{{g}^{2}}-{{z}^{2}}} &\text{ if }& \phi \le 0\,.
\end{array} \right.
\EE
Moreover, $c_1-c_2=(c_1-c_2)(\Sigma(\phi))$ is also a function of $\phi$. Note that $\Sigma(0)={{\Sigma}_{z}}$, $\Sigma^{\prime}(0)=0$ and $\left( {{c}_{1}}-{{c}_{2}} \right)\left( \Sigma \left( 0 \right) \right)=\left( {{c}_{1}}-{{c}_{2}} \right)\left( {{\Sigma }_{z}} \right)=0$. Then (\ref{cn-cp}) and (\ref{drho-phi}) imply
	 \[\frac{d}{d\phi }\left( {{c}_{1}}-{{c}_{2}} \right)=\frac{d}{d\Sigma }\left( {{c}_{1}}-{{c}_{2}} \right)\frac{d\Sigma }{d\phi }=q\frac{\Sigma {{e}^{\left( g+z \right)\Sigma }}+2\left( g+z \right)}{\left( 1+g\Sigma \right){{e}^{\left( g+z \right)\Sigma }}+{{g}^{2}}-{{z}^{2}}}>0\quad \text{for}\;\phi \in \mathbb{R}.\]
Therefore, $c_1-c_2$ is strictly monotone increasing to $\phi$ and we complete the proof.
\end{proof}

\begin{figure}[htp]
	\centering
	\includegraphics[scale=0.4, bb=0 0 650 300]{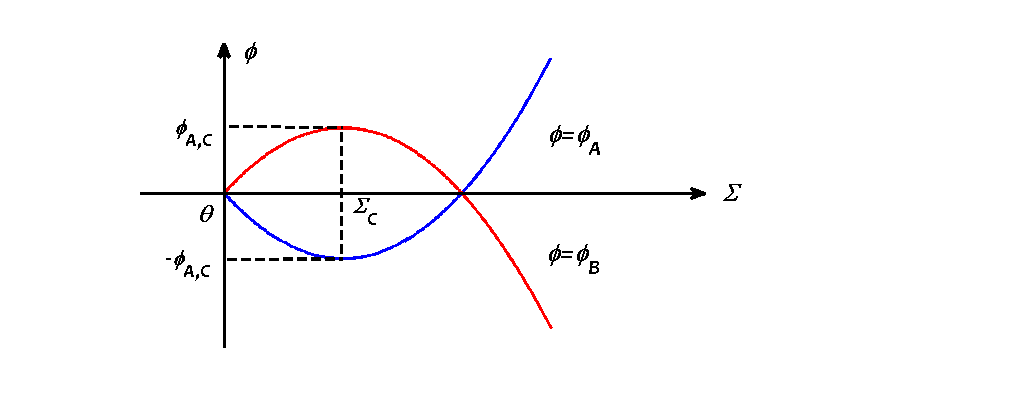}
	\vspace*{-1cm}
  	\caption{$\theta = (\Sigma_z,0)$ in $(\Sigma,\phi)$ coordinates}
 	\label{fig2}
\end{figure}

When $z=g_{12}$ is increased, for example when the ion is divalent like calcium, the profiles of $\phi_A$ and $\phi_B$ may lose monotonicity and become oscillatory. It is well known in experiments that calcium has profound and complex effects on the current voltage relations of channels (cf.~\cite{AM-pnas90, FH-jp57}). Suppose $z>\sqrt{1+{{g}^{2}}}>0$. Then $z^2-g^2>1$ and there exists a unique ${{\Sigma}_{c}}>0$ (because  $\left( 1+g\Sigma \right){{e}^{\left( g+z \right)\Sigma }}$ is strictly monotone increasing to $\Sigma>0$) depending on $\Sigma_z$ such that
\[(1+g{{\Sigma }_{c}}){{e}^{(g+z){{\Sigma }_{c}}}}+{{g}^{2}}-{{z}^{2}}=0.\]
Note that $\frac{d{{\phi }_{A}}}{d\Sigma }\left( {{\Sigma }_{c}} \right)=\frac{d{{\phi }_{B}}}{d\Sigma }\left( {{\Sigma }_{c}} \right)=0$ if $\Sigma_c>\Sigma_z>0$. We shall prove that $\Sigma_c$ may be located in the domain of $\phi_A$ and $\phi_B$ i.e. $\Sigma_c>\Sigma_z>0$ if $z$ is sufficiently large (see Proposition~\ref{rho_c}). By (\ref{dphia}) and (\ref{dphib}), $\frac{d{{\phi }_{A}}}{d\Sigma }<0$ on $\left( {{\Sigma }_{z}},{{\Sigma }_{c}} \right)$, $\frac{d{{\phi }_{A}}}{d\Sigma }>0$ on $\left( {{\Sigma }_{c}},\infty\right)$, $\frac{d{{\phi }_{B}}}{d\Sigma }>0$ on $\left( {{\Sigma }_{z}},{{\Sigma }_{c}} \right)$, $\frac{d{{\phi }_{B}}}{d\Sigma }<0$ on $\left( {{\Sigma }_{c}},\infty\right)$. Then $\Sigma_c$ is a unique (global) minimal point of $\phi_A$ and a unique (global) maximal point of $\phi_B$, respectively (see Figure~\ref{fig2}). Moreover, by~(\ref{sym}),
\BE\label{pi-ac1}
{{\phi }_{A,c}}\equiv -{{\phi }_{A}}\left( {{\Sigma }_{c}} \right)=-\underset{\Sigma >{{\Sigma }_{z}}}{\mathop{\min }}\,{{\phi }_{A}}\left( \Sigma  \right)=\underset{\Sigma >{{\Sigma }_{z}}}{\mathop{\max }}\,{{\phi }_{B}}\left( \Sigma  \right)={{\phi }_{B}}\left( {{\Sigma }_{c}} \right)>0\,.
\EE
\vspace*{1.2cm}
\vspace*{-0.3cm}

\begin{figure}[htp]
 	\centering
 	\includegraphics[scale=0.4, bb=0 0 1000 350]{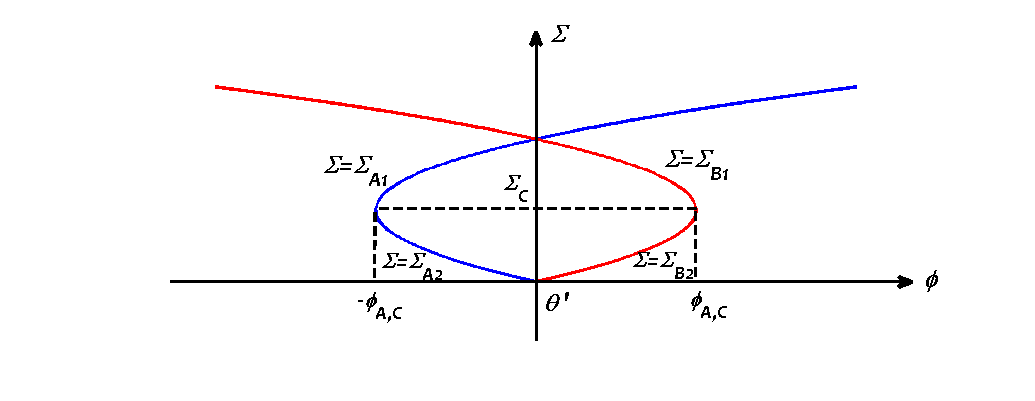}
	\vspace*{-1cm}  	
  	\caption{$\theta' =(0,\Sigma_z)$}
  	\label{fig3}
\end{figure}

\noindent By Figure~\ref{fig2}, the inverse image of function ${{\phi }_{A}}$ consists of two functions ${{\Sigma }_{A_1}}:\left( -{{\phi }_{A,c}},\infty \right)\to \left( {{\Sigma }_{c}},\infty  \right)$ and ${{\Sigma }_{{{A}_{2}}}}:\left[ -{{\phi }_{A,c}},0 \right]\to \left[ {{\Sigma }_{z}},{{\Sigma }_{c}} \right]$ such that $\frac{d{{\Sigma}_{A_1}}}{d\phi }>0$ on $\left( -{{\phi }_{A,c}},\infty \right)$ and $\frac{d{{\Sigma}_{A_2}}}{d\phi }<0$ on $\left( -{{\phi}_{A,c}},0 \right)$ (see Figure~\ref{fig3}).

\vspace*{0.5cm}

\noindent Moreover, by (\ref{dphia}),
\BE\label{d-rho-a+}
\frac{d{{\Sigma }_{{{A}_{1}}}}}{d\phi }=q\,\frac{{{e}^{(g+z){{\Sigma }_{{{A}_{1}}}}}}\sqrt{\Sigma _{{{A}_{1}}}^{2}-4{{e}^{-(g+z){{\Sigma }_{{{A}_{1}}}}}}}}{(1+g{{\Sigma }_{{{A}_{1}}}}){{e}^{(g+z){{\Sigma }_{{{A}_{1}}}}}}+{{g}^{2}}-{{z}^{2}}}>0\quad \text{ for }\quad \phi >-{{\phi }_{A,c}}\,,
\EE
and
$$
\frac{d{{\Sigma }_{{{A}_{2}}}}}{d\phi }=q\,\frac{{{e}^{(g+z){{\Sigma }_{{{A}_{2}}}}}}\sqrt{\Sigma _{{{A}_{2}}}^{2}-4{{e}^{-(g+z){{\Sigma }_{{{A}_{2}}}}}}}}{(1+g{{\Sigma }_{{{A}_{2}}}}){{e}^{(g+z){{\Sigma }_{{{A}_{2}}}}}}+{{g}^{2}}-{{z}^{2}}}<0\quad \text{ for }\quad -{{\phi }_{A,c}}<\phi <0\,.
$$
Similarly, the inverse image of function ${{\phi }_{B}}$ consists of another two functions ${{\Sigma }_{B_1}}:\left( -\infty ,{{\phi }_{A,c}} \right)\to \left( {{\Sigma }_{c}},\infty\right)$ and  ${{\Sigma }_{{{B}_{2}}}}:\left[ 0,{{\phi }_{A,c}} \right]\to \left[ {{\Sigma }_{z}},{{\Sigma }_{c}} \right]$ such that $\frac{d{{\Sigma }_{B_1}}}{d\phi }<0$ on   $\left( -\infty ,{{\phi }_{A,c}} \right)$ and $\frac{d{{\Sigma }_{B_2}}}{d\phi }>0$ on $\left( 0,{{\phi }_{A,c}} \right)$. Moreover, by (\ref{dphib}),
\BE\label{d-rho-b+}
\frac{d{{\Sigma }_{{{B}_{1}}}}}{d\phi }=-q\,\frac{{{e}^{(g+z){{\Sigma }_{{{B}_{1}}}}}}\sqrt{\Sigma _{{{B}_{1}}}^{2}-4{{e}^{-(g+z){{\Sigma }_{{{B}_{1}}}}}}}}{(1+g{{\Sigma }_{{{B}_{1}}}}){{e}^{(g+z){{\Sigma }_{{{B}_{1}}}}}}+{{g}^{2}}-{{z}^{2}}}<0\quad \text{ for }\quad \phi <{{\phi }_{A,c}}\,,
\EE
and
\[\frac{d{{\Sigma }_{{{B}_{2}}}}}{d\phi }=-q\,\frac{{{e}^{(g+z){{\Sigma }_{{{B}_{2}}}}}}\sqrt{\Sigma _{{{B}_{2}}}^{2}-4{{e}^{-(g+z){{\Sigma }_{{{B}_{2}}}}}}}}{(1+g{{\Sigma }_{{{B}_{2}}}}){{e}^{(g+z){{\Sigma }_{{{B}_{2}}}}}}+{{g}^{2}}-{{z}^{2}}}>0\quad \text{ for }\quad 0<\phi <{{\phi }_{A,c}}\,.\]
Thus by~(\ref{cn-cp}), we may consider two functions of  $({{c}_{1}}-{{c}_{2}})\circ {{\Sigma }_{{{A}_{1}}}}$ and  $({{c}_{1}}-{{c}_{2}})\circ {{\Sigma }_{B_1}}$ as follows:
\BE\label{cn-cp-a+}
({{c}_{1}}-{{c}_{2}})({{\Sigma }_{{{A}_{1}}}}(\phi ))=\sqrt{\Sigma_{{{A}_{1}}}^{2}-4{{e}^{-(g+z){{\Sigma }_{{{A}_{1}}}}}}}\quad \hbox{for}\quad \phi \ge -{{\phi }_{A,c}}\,,
\EE
and
\BE\label{cn-cp-b+}
(c_1-c_2)(\Sigma_{B_1}(\phi))=-\sqrt{\Sigma_{B_1}^2-4e^{-(g+z)\Sigma_{B_1}}}
\quad\hbox{ for }\quad\phi\leq \phi_{A,c}\,.
\EE
Note that $(c_1-c_2)(\Sigma_{A_1}(\cdot))$ and $(c_1-c_2)(\Sigma_{B_1}(\cdot))$ are continuous functions on $[-\phi_{A,c},\phi_{A,c}]$. Moreover, by~(\ref{d-rho-a+})-(\ref{cn-cp-b+}), we have
\BE\label{d-cn-cp-a+}
\frac{d}{d\phi }({{c}_{1}}-{{c}_{2}})({{\Sigma }_{{{A}_{1}}}}(\phi ))=q\,\frac{{{e}^{(g+z){{\Sigma }_{{{A}_{1}}}}}}[{{\Sigma }_{{{A}_{1}}}}+2(g+z){{e}^{-(g+z){{\Sigma }_{{{A}_{1}}}}}}]}{(1+g{{\Sigma }_{{{A}_{1}}}}){{e}^{(g+z){{\Sigma }_{{{A}_{1}}}}}}+{{g}^{2}}-{{z}^{2}}}>0\quad \text{ for }\quad \phi >-{{\phi }_{A,c}}\,,
\EE
and
\BE\label{d-cn-cp-b+}
\frac{d}{d\phi }({{c}_{1}}-{{c}_{2}})({{\Sigma }_{{{B}_{1}}}}(\phi ))=q\,\frac{{{e}^{(g+z){{\Sigma }_{{{B}_{1}}}}}}[{{\Sigma }_{{{B}_{1}}}}+2(g+z){{e}^{-(g+z){{\Sigma }_{{{B}_{1}}}}}}]}{(1+g{{\Sigma }_{{{B}_{1}}}}){{e}^{(g+z){{\Sigma }_{{{B}_{1}}}}}}+{{g}^{2}}-{{z}^{2}}}>0\quad \text{ for }\quad \phi <{{\phi }_{A,c}}\,.
\EE
Here we have used~(\ref{rho_c2}) and~(\ref{rho_c3}). Consequently, $(c_1-c_2)(\Sigma_{A_1}(\cdot))$ and $(c_1-c_2)(\Sigma_{B_1}(\cdot))$ are smooth functions on $(-\phi_{A,c},\phi_{A,c})$. Since $(c_1-c_2)(\Sigma_{A_1}(\cdot))$ and $(c_1-c_2)(\Sigma_{B_1}(\cdot))$ are strictly monotone increasing to $\phi$ (see~(\ref{d-cn-cp-a+}) and (\ref{d-cn-cp-b+})), then we may use~(\ref{cn-cp}) to get
\BE\label{v-cn-cp-a+}
({{c}_{1}}-{{c}_{2}})({{\Sigma }_{{{A}_{1}}}}(\phi ))\ge ({{c}_{1}}-{{c}_{2}})({{\Sigma }_{{{A}_{1}}}}(-{{\phi }_{A,c}}))=\sqrt{\Sigma _{c}^{2}-4{{e}^{-(g+z){{\Sigma }_{c}}}}}>0\,,
\EE
\BE\label{v-cn-cp-b+}
({{c}_{1}}-{{c}_{2}})({{\Sigma }_{{{B}_{1}}}}(\phi ))\le ({{c}_{1}}-{{c}_{2}})({{\Sigma }_{{{B}_{1}}}}({{\phi }_{A,c}}))=-\sqrt{\Sigma_{c}^{2}-4{{e}^{-(g+z){{\Sigma }_{c}}}}}<0\,,
\EE
for $\phi\in (-\phi_{A,c},\phi_{A,c})$.

Now we claim that if $z$ is sufficiently large, then $\Sigma_c>\Sigma_z>0$ i.e. $\Sigma_c$ is located in the domain of $\phi_A$ and $\phi_B$ as follows:
\medskip
\noindent

\BP\label{rho_c} Let \BE\label{gc1} g_c=\inf\{z>\sqrt{1+g^2}:
\hbox{ there\ exists}\: \Sigma_{c,z}>\Sigma_z>0\:\hbox{ such\ that}\:
(1+g\Sigma_{c,z})e^{(g+z)\Sigma_{c,z}}+g^2-z^2=0 \}\,,\EE where $\Sigma_z>0$ is
the unique solution of $\Sigma=2e^{-\frac{1}{2}(g+z)\Sigma}$ for
$z>0$. Then for $z>g_c$, there exists a unique $\Sigma_c=\Sigma_{c,z}>\Sigma_z$ depending on $z$ such that
$(1+g\Sigma_c)e^{(g+z)\Sigma_c}+g^2-z^2=0$. Conversely, for
$0<z<g_c$, no such $\Sigma_c$ exists and
$(1+g\Sigma)e^{(g+z)\Sigma}+g^2-z^2>0$ for
$\Sigma\geq\Sigma_z>0$.
\EP
\begin{proof}
Firstly, we claim that $g_c$ is well-defined. For any $z>0$, we
may define a function $f_z=f_z(\Sigma)$ by \BE\label{fun_z}
f_z(\Sigma)=(1+g\Sigma)e^{(g+z)\Sigma}+g^2-z^2\quad\hbox{ for }\quad
\Sigma>0\,.\EE Then it is obvious that $f_z(+\infty)=\infty$,
\BE\label{mf}
f_z'(\Sigma)=[g+(1+g\Sigma)(g+z)]e^{(g+z)\Sigma}>0\quad\hbox{ for }\quad
\Sigma, z>0\,,\EE and $f_z(0)=1+g^2-z^2<0$ if $z>\sqrt{1+g^2}$.
Hence there exists a unique $\Sigma_{c,z}>0$ such that
$f_z(\Sigma_{c,z})=0$. Let $\Sigma///_z>0$ be the unique solution of
\BE\label{rho_z} \Sigma_z=2e^{-\frac{1}{2}(g+z)\Sigma_z}\quad\hbox{
for }\quad z>0\,.\EE

Now we prove $\Sigma_{c,z}>\Sigma_z$ as $z$ sufficiently large. By~(\ref{rho_z}), $\Sigma_z$ is decreasing to $z$ (differentiate (\ref{rho_z}) to $z$) and $z=-\left(g+\frac{2\ln\Sigma_z-\ln 4}{\Sigma_z}\right)$. Thus $\Sigma_z\to 0$ as $z\to\infty$ and
\begin{eqnarray*} f_z(\Sigma_z)&=&(1+g\Sigma_z)e^{(g+z)\Sigma_z}+g^2-z^2 \\
&=&[4(1+g\Sigma_z)+(g^2-z^2)\Sigma_z^2]/\Sigma_z^2 \quad\hbox{by~(\ref{rho_z})}\\
&=&[4(1+g\Sigma_z)-2g\Sigma_z(2\ln\Sigma_z-\ln 4)-(2\ln\Sigma_z-\ln
4)^2]/\Sigma_z^2 \to -\infty\quad\hbox{ as }\:z\to\infty\,,
\end{eqnarray*}
and then $f_z(\Sigma_z)<0$ as $z$ sufficiently large. Since $f_z(\Sigma_{c,z})=0$ and $f_z(\Sigma_z)<0$ as $z$ sufficiently large, then by~(\ref{mf}), we have $\Sigma_{c,z}>\Sigma_z$ as $z$ sufficiently large. Consequently, the set
\begin{align}\label{gc2}\mathcal{Z}&=\{z>\sqrt{1+g^2}: \exists
\Sigma_{c,z}>\Sigma_z>0\:\hbox{ such\ that}\: f_z(\Sigma_{c,z})=0 \} \\
&=\{z>\sqrt{1+g^2}: f_z(\Sigma_z)<0\} \notag
\end{align}
is nonempty and the value ${{g}_{c}}=\underset{z\in \mathcal{Z}}{\mathop{\inf }}\,z$ (defined in (\ref{gc1})) is well-defined. Note that the existence of $\Sigma_{c,z}$ with $f_z(\Sigma_{c,z})=0$ is guaranteed due to $z>\sqrt{1+{{g}^{2}}}$, so (\ref{mf}) implies $\Sigma_{c,z}>\Sigma_z$ if $f_z(\Sigma_z)<0$ holds true.

To complete the proof of Proposition~\ref{rho_c}, we need the following result: \\
{\bf Claim~1.}~~{\it Suppose $f_{z_0}(\Sigma_{z_0})=0$ and $\Sigma_{z_0}>0$ for some $z_0>\sqrt{1+g^2}$. Then
there exist $z_l, z_r>\sqrt{1+g^2}$ and $z_l<z_0<z_r$ such that $f_{z}(\Sigma_{z})>0$ for $z\in (z_l,z_0)$ and $f_{z}(\Sigma_{z})<0$ for $z\in (z_0,z_r)$.}
\begin{proof}
By~(\ref{fun_z}) and (\ref{rho_z}),
\BE\label{2.30-1}
f\left( {{\Sigma }_{z}} \right)=\frac{4\left( 1+g{{\Sigma }_{z}} \right)}{\Sigma _{z}^{2}}+{{g}^{2}}-{{z}^{2}}\,.
\EE
Then $f_{z_0}(\Sigma_{z_0})=0$ gives
\BE\notag
4\frac{1+g\Sigma_{z_0}}{\Sigma_{z_0}^2}= z_0^2-g^2\,,
\EE
and $\Sigma_{z_0}$ satisfies $(z_0^2-g^2)\Sigma_{z_0}^2-4g\Sigma_{z_0}-4=0$ having solutions as $\Sigma_{z_0}=\frac{2}{z_0-g}$ and $\Sigma_{z_0}=-\frac{2}{z_0+g}$. Hence due to $\Sigma_{z_0}>0$,
\BE\label{rho_z0-1}
\Sigma_{z_0}=\frac{2}{z_0-g}\,.
\EE
Note that $z_0>\sqrt{1+g^2}>\pm g$. Differentiating (\ref{rho_z}) and (\ref{2.30-1}) to $z$, we have
\begin{align}\nonumber
\frac{d}{dz}f_z(\Sigma_z)=&
-4\frac{2+g\Sigma_z}{\Sigma_z^3}\,\frac{d\Sigma_z}{dz}-2z\,, \\
&\notag \\
\notag \frac{d\Sigma_z}{dz}=& \frac{-\Sigma_z^2}{(g+z)\Sigma_z+2}\,.
\end{align}
Thus by~(\ref{rho_z0-1}), we obtain
\BE\label{f'z_0-1} \frac{d}{dz}f_z(\Sigma_z)|_{z=z_0}=-z_0-g<0\,.\EE
Therefore, by~(\ref{f'z_0-1}), we may complete the proof of Claim~1.
\end{proof}

It is obvious that
\BE\label{z-small} f_z(\Sigma)>0 \quad\hbox{for}\quad
\Sigma>0\quad\hbox{ and }\quad 0<z\leq\sqrt{1+g^2}\,.\EE
Now we want to prove that \BE\label{setZ}
\mathcal{Z}=(g_c,\infty)\,, \EE where ${{g}_{c}}=\underset{z\in \mathcal{Z}}{\mathop{\inf }}\,z$. Due to the continuity of $f_z$, (\ref{gc2}) implies that the set $\mathcal{Z}$ is open. Suppose the set $\mathcal{Z}$ has two components. Then without loss of generality, we may assume that there exists $z_a>g_c$ such that $\mathcal{Z}=(g_c,z_a)\cup (z_a,\infty)$. Hence $f_{z_a}(\Sigma_{z_a})=0$ and $f_z(\Sigma_z)<0$ for $z\in (g_c,z_a)\cup (z_a,\infty)$. However, Claim~1 implies that $f_{z}(\Sigma_{z})>0$ for $z\in (z_l,z_a)$ which contradicts to $f_z(\Sigma_z)<0$ for $z\in (g_c,z_a)$. Thus~the proof of (\ref{setZ}) is done. On the other hand, Claim~1 also implies that
\BE\label{z-small-1} f_z(\Sigma_z)>0\quad\hbox{ for }\quad 0<z<g_c\,.\EE
Otherwise, by~(\ref{z-small}), there exists $z_b\in (\sqrt{1+g^2},g_c)$ such that $f_{z_b}(\Sigma_{z_b})=0$. Then as for~(\ref{f'z_0-1}), we have $\frac{d}{dz}f_z(\Sigma_z)|_{z=z_b}=-z_b-g<0$ and hence there exists $z_c\in (z_b, g_c)$ such that $f_{z_c}(\Sigma_{z_c})<0$ which contradicts to~(\ref{setZ}). Therefore, by~(\ref{mf}) and (\ref{z-small-1}), we complete the proof of Proposition~\ref{rho_c}.
\end{proof}

\BR\label{gc}
\medskip
\noindent

\begin{enumerate} \item[(i)]~~The proof of Proposition~\ref{rho_c} shows that $f_z(\Sigma_z)>0$ for $0<z<g_c$
and $f_z(\Sigma_z)<0$ for $z>g_c$ (see~(\ref{setZ}) and
(\ref{z-small-1})). Hence by the continuity of $f_z$,
$f_{g_c}(\Sigma_{g_c})=0$. \item[(ii)]~~By (\ref{mf}) and
(\ref{z-small-1}), we have \BE\label{mf1}
f_z(\Sigma)=(1+g\Sigma)e^{(g+z)\Sigma}+g^2-z^2>0\quad\hbox{for}\quad
\Sigma\geq\Sigma_z\quad\hbox{ and }\quad 0<z<g_c\,. \EE
\item[(iii)]~~By (\ref{fun_z}), ${{f}_{z}}\left( {{\Sigma }_{z}} \right)>0$ as $z=\sqrt{1+{{g}^{2}}}$ but $f_{g_c}(\Sigma_{g_c})=0$. Hence Remark~\ref{gc}~(i) implies $g_c>\sqrt{1+g^2}$.
\end{enumerate} \ER

Suppose $0<z<g_c$. Then (\ref{mf1}) gives $f_{z}(\Sigma)>0$ for $\Sigma\geq\Sigma_{z}$. Hence
by~(\ref{dphia}) and (\ref{dphib}), $\frac{d\phi_A}{d\Sigma}>0$ and
$\frac{d\phi_B}{d\Sigma}<0$ for $\Sigma\geq\Sigma_{z}$ which gives  ${{\phi }_{A}}(\Sigma )>{{\phi }_{A}}({{\Sigma }_{z}})=0={{\phi }_{B}}({{\Sigma }_{z}})>{{\phi }_{B}}(\Sigma )$ for $\Sigma>\Sigma_{z}$. Thus as for Proposition~\ref{mon}, $\Sigma=\Sigma(\phi)$ can be a single-valued function of $\phi$ with domain as the entire space $\mathbb{R}$ and range  $\left[ {{\Sigma }_{z}},\infty  \right)$ such that $\Sigma(0)=\Sigma_{z}$ and $c_1-c_2$ is strictly monotone increasing to $\phi$. Moreover, $\Sigma\to\infty$ as $\phi\to\pm\infty$ and $(c_1-c_2)(\Sigma(\phi))\to\pm\infty$ as $\phi\to\pm\infty$.

Suppose $z>g_c>0$. Then Proposition~\ref{rho_c} gives that there exists a unique  ${{\Sigma }_{c}}\in ({{\Sigma }_{z}},\infty )$ such that
\BE\label{rho_c1}
(1+g{{\Sigma }_{c}}){{e}^{(g+z){{\Sigma }_{c}}}}+{{g}^{2}}-{{z}^{2}}=0\,,\EE
which implies
\BE \label{rho_c2}
(1+g\Sigma ){{e}^{(g+z)\Sigma }}+{{g}^{2}}-{{z}^{2}}>0,\quad \hbox{for}\quad \Sigma >{{\Sigma }_{c}}\,,\EE
and
\BE \label{rho_c3}
(1+g\Sigma ){{e}^{(g+z)\Sigma }}+{{g}^{2}}-{{z}^{2}}<0,\quad \hbox{for}\quad {{\Sigma }_{z}}\le \Sigma <{{\Sigma }_{c}}\,,
\EE By (\ref{ph1}) and (\ref{dphia}), we have $\frac{d\phi_A}{d\Sigma}>0$ for $\Sigma>\Sigma_c$; $\frac{d\phi_A}{d\Sigma}<0$ for $\Sigma_{z}<\Sigma<\Sigma_c$, and $\phi_A$ tends to $+\infty$ as $\Sigma$ goes to $+\infty$. Hence $\Sigma_c$ is the unique minimum point of $\phi_A$. Since $\Sigma_{z}^{2}=4{{e}^{-(g+z){{\Sigma }_{z}}}}$, then  ${{\phi }_{A}}({{\Sigma }_{z}})=0$ which implies $-\phi_{A,c}=\phi_A(\Sigma_c)<0$. Since $\Sigma_c$ satisfies  $(1+g{{\Sigma }_{c}}){{e}^{(g+z){{\Sigma }_{c}}}}={{z}^{2}}-{{g}^{2}}$ i.e.~  ${{\Sigma }_{c}}+\frac{\ln \left( 1+g{{\Sigma }_{c}} \right)}{g+z}=\frac{\ln \left( {{z}^{2}}-{{g}^{2}} \right)}{g+z}$, then $\Sigma_c$ must tend to zero as $z$ goes to infinity. Note that $g>0$ is a fixed constant. Consequently,  $-\ln \left[ \frac{1}{2}\left( {{\Sigma }_{c}}+\sqrt{\Sigma_{c}^{2}-4{{e}^{-(g+z){{\Sigma }_{c}}}}} \right) \right]\to +\infty $ as $z\to +\infty$, and then
\begin{align}
q\,\phi_{A,c}&=q\,\phi_A(\Sigma_c) \notag \\
&=-\ln \left[ \frac{1}{2}\left( {{\Sigma }_{c}}+\sqrt{\Sigma_{c}^{2}-4{{e}^{-(g+z){{\Sigma }_{c}}}}} \right) \right]+\frac{g+z}{2}{{\Sigma }_{c}}+\frac{g-z}{2}\sqrt{\Sigma_{c}^{2}-4{{e}^{-(g+z){{\Sigma }_{c}}}}}\notag
\\
&=-\ln \left[ \frac{1}{2}\left( {{\Sigma }_{c}}+\sqrt{\Sigma_{c}^{2}-4{{e}^{-(g+z){{\Sigma }_{c}}}}} \right) \right]+\frac{g}{2}\left( {{\Sigma }_{c}}+\sqrt{\Sigma_{c}^{2}-4{{e}^{-(g+z){{\Sigma }_{c}}}}} \right)
\notag
\\
&~~~~~+\frac{z}{2}\left( {{\Sigma }_{c}}-\sqrt{\Sigma_{c}^{2}-4{{e}^{-(g+z){{\Sigma }_{c}}}}} \right)\notag
\\
&\geq -\ln \left[ \frac{1}{2}\left( {{\Sigma }_{c}}+\sqrt{\Sigma_{c}^{2}-4{{e}^{-(g+z){{\Sigma }_{c}}}}} \right) \right]\to +\infty
\notag \end{align}
as $z\to +\infty$.  Thus $\phi_{A,c}\to +\infty$ as $z\to +\infty$ and $g>0$ is fixed. Besides, since   ${{e}^{-(g+z){{\Sigma }_{c}}}}=(1+g{{\Sigma }_{c}})/({{z}^{2}}-{{g}^{2}})$ and $\Sigma_c\to 0$ as $z\to\infty$, then by~(\ref{cn-cp}), we have $(c_1-c_2)(\Sigma_c)\to 0$ as $z\to +\infty$ and $g>0$ is fixed. Therefore, we may summarize the above results as follows: \BT\label{pro-cn-cp}
\medskip
\noindent
\begin{enumerate} \item[(i)]~~Suppose $0<z<g_c$. Then $(c_1-c_2)\circ\Sigma$ is a monotone increasing function to $\phi\in \mathbb{R}$ satisfying $(c_1-c_2)(\Sigma(\phi))\to\pm\infty$ as $\phi\to\pm\infty$, respectively. \item[(ii)]~~Suppose $z>g_c$. Then there are two functions $\Sigma_{A_1}$ and $\Sigma_{B_1}$ such that $(c_1-c_2)\circ\Sigma_{A_1}:[-\phi_{A,c},\infty)\to\mathbb{R}$ and $(c_1-c_2)\circ\Sigma_{B_1}:(-\infty,\phi_{A,c}]\to\mathbb{R}$ are monotone increasing functions of $\phi$, where $\phi_{A,c}$ satisfies $\phi_{A,c}\to +\infty$ and $(c_1-c_2)(\Sigma_c)\to 0$ as $z\to +\infty$ and $g>0$ is fixed. Moreover, \begin{align}\notag & (c_1-c_2)\circ\Sigma_{A_1}(-\phi_{A,c})=(c_1-c_2)(\Sigma_c)>0\,,\\
&
(c_1-c_2)\circ\Sigma_{B_1}(\phi_{A,c})=-(c_1-c_2)(\Sigma_c)<0\,,\notag\end{align}
$$\lim_{\phi\to\infty}(c_1-c_2)\circ\Sigma_{A_1}(\phi)=\infty\quad\hbox{and}\: \lim_{\phi\to-\infty}(c_1-c_2)\circ\Sigma_{B_1}(\phi)=-\infty\,.$$
\end{enumerate} Here $\circ$ denotes function composition and $g_c$ is the positive constant defined in Proposition~\ref{rho_c}. \ET

\section{Proof of Theorem~\ref{thm1.1} and~\ref{thm1.2} }\label{mul-solu-3s}

\subsection{ Proof of Theorem~\ref{thm1.1} }\label{mul-solu-3s-1}
In this section, we study multiple solutions of the system of equations (\ref{eqn3-4})-(\ref{bdc-Rbn-0}) with $N=3$ and the following assumptions:
\BE\label{eqn-3.0}
\rho_0>0,\quad g_{i3}=g_{3i}=0\quad\hbox{for}\quad i=1,2,3\,.
\EE
Then we may get solutions of (\ref{eqn3-4}) by solving
\BE\label{eqn-3.1}
\ln {{c}_{i}}+{{z}_{i}}\phi +\sum\limits_{j=1}^{2}{{{g}_{ij}}{{c}_{j}}=0}\quad\hbox{for}\quad i=1,2\,,
\EE
and let
\BE\label{eqn-3.2}
{{c}_{3}}={{e}^{-{{z}_{3}}\phi }}\,.
\EE
Note that (\ref{eqn-3.1}) is same as (\ref{eqn3-4}) with $N=2$. Assume
\BE\label{eqn-3.2.1}
z_2=-z_1=q\geq 1\,,\quad g_{11}=g_{22}>0\quad\hbox{ and }\quad g_{12}>g_c>0\,,
\EE
where $g_c>0$ is a sufficiently large constant defined in Proposition~\ref{rho_c}. We shall use (\ref{eqn-3.2.1}) and set $\rho_0>0$ in order to apply Theorem~\ref{pro-cn-cp}~(ii) (in Section~\ref{ss-solu-2s}) and Lemma~\ref{bd-solu} (in Section~\ref{sec4}) for the proof of Theorem~\ref{thm1.1} which gives multiple solutions of (\ref{eqn3-4})-(\ref{bdc-Rbn-0}) with $N=3$ and $\rho_0>0$.

\begin{figure}[h]
 	\centering
 	\includegraphics[scale=0.25, bb=0 0 800 307]{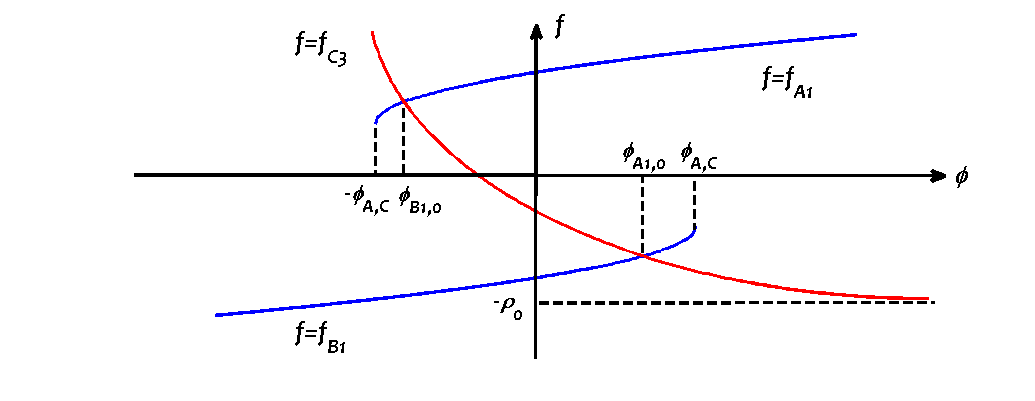}
	\vspace*{-0.6cm}  	
  	\caption{Figures of $f_{A_1}$, $f_{B_1}$ and $f_{c_3}$}
  	\label{fig4}
\end{figure}

By Theorem~\ref{pro-cn-cp}~(ii), equation (\ref{eqn-3.1}) has multiple solutions
\BE\label{multi-sol1}
\left( {{c}_{1}},{{c}_{2}} \right)=\left( {{c}_{1}}\left( {{\Sigma }_{{{A}_{1}}}}\left( \phi  \right) \right),{{c}_{2}}\left( {{\Sigma }_{{{A}_{1}}}}\left( \phi\right) \right) \right)
\quad\hbox{and}\quad
\left( {{c}_{1}},{{c}_{2}} \right)=\left( {{c}_{1}}\left( {{\Sigma }_{{{B}_{1}}}}\left( \phi  \right) \right),{{c}_{2}}\left( {{\Sigma }_{{{B}_{1}}}}\left( \phi\right) \right) \right)
\EE
such that ${{f}_{{{A}_{1}}}}(\phi )=q\,({{c}_{1}}-{{c}_{2}})\left( {{\Sigma }_{{{A}_{1}}}}(\phi ) \right)$ and
${{f}_{{{B}_{1}}}}(\phi )=q\,({{c}_{1}}-{{c}_{2}})\left( {{\Sigma }_{{{B}_{1}}}}(\phi ) \right)$ are monotone increasing to $\phi$ but the values of $f_{A_1}$ and $f_{B_1}$ are away from zero (see Figure~\ref{fig4}). By Lemma~\ref{ub-solu}, it is impossible to get uniformly bounded solution by solving either $\varepsilon\phi^{\prime\prime}(x)=f_{A_1}(\phi(x))$ or $\varepsilon\phi^{\prime\prime}(x)=f_{B_1}(\phi(x))$ for $x\in (-1,1)$. This motivates us to develop Lemma~\ref{bd-solu} (in Section~\ref{sec4}), and use (\ref{eqn-3.2}) to transform (\ref{eqn-1.9}) into the following equations:
\BE\label{eqn3-5-A} \varepsilon\phi^{\prime\prime}(x)=f_A(\phi(x)) \quad\hbox{for }\quad x\in (-1,1)\,, \EE
and
\BE\label{eqn3-5-B}
\varepsilon\phi^{\prime\prime}(x)=f_B(\phi(x)) \quad\hbox{for
}\quad x\in (-1,1)\,, \EE
where
\[{{f}_{A}}(\phi )=q\,({{c}_{1}}-{{c}_{2}})\left( {{\Sigma }_{{{A}_{1}}}}(\phi ) \right)-{{z}_{3}}{{e}^{-{{z}_{3}}\phi }}+{{\rho }_{0}}\,,\]
and	
\[{{f}_{B}}(\phi )=q\,({{c}_{1}}-{{c}_{2}})\left( {{\Sigma }_{{{B}_{1}}}}(\phi ) \right)-{{z}_{3}}{{e}^{-{{z}_{3}}\phi }}+{{\rho }_{0}}\,.\]
We may denote $f_A$ and $f_B$ as follows:
${{f}_{A}}\left(\phi
\right)={{f}_{{{A}_{1}}}}\left(\phi\right)-{{f}_{{{c}_{3}}}}\left(\phi\right)$ and
${{f}_{B}}\left(\phi\right)={{f}_{{{B}_{1}}}}\left(\phi\right)-{{f}_{{{c}_{3}}}}\left(\phi  \right)$,  where ${{f}_{{{A}_{1}}}}(\phi )=q\,({{c}_{1}}-{{c}_{2}})\left( {{\Sigma }_{{{A}_{1}}}}(\phi ) \right)$,
${{f}_{{{B}_{1}}}}(\phi )=q\,({{c}_{1}}-{{c}_{2}})\left( {{\Sigma }_{{{B}_{1}}}}(\phi ) \right)$, and  ${{f}_{{{c}_{3}}}}(\phi )={{z}_{3}}{{e}^{-{{z}_{3}}\phi }}-{{\rho }_{0}}$.

Let $\rho_0>0$. Then Theorem~\ref{pro-cn-cp}~(ii) (in Section~\ref{ss-solu-2s}) implies that as $g_{12}=z\geq g_{\rho_0}>g_c>0$ ($g_{\rho_0}$ is a large constant depending on $\rho_0$), both functions $f_{A_1}$ and $f_{B_1}$ intersect with the function $f_{c_3}$ at $\phi_{A_1,0}$ and $\phi_{B_1,0}$, respectively (see Figure~\ref{fig4}). Note that the assumption $\rho_0>0$ is necessary for the existence of $\phi_{A_1,0}$ and $\phi_{B_1,0}$.
%\vspace*{0.1cm}
Moreover, $f_A=f_{A_1}-f_{c_3}$ and $f_B=f_{B_1}-f_{c_3}$ satisfy
\begin{enumerate}
\item[(1)]~~$f_A:[-\phi_{A,c},\infty)\to\mathbb{R}$ is smooth and strictly monotone increasing, $-\phi_{A,c}<0, f_A(-\phi_{A,c})<0, f_A(\infty)>0$ and $f_A(\phi_{{A_1},0})=0$ for some $\phi_{{A_1},0}>-\phi_{A,c}$.
\item[(2)]~~$f_B:(-\infty,\phi_{A,c}]\to\mathbb{R}$ is smooth and strictly monotone increasing, $\phi_{A,c}>0, f_B(\phi_{A,c})>0, f_B(-\infty)<0$ and $f_B(\phi_{{B_1},0})=0$ for some $\phi_{{B_1},0}<\phi_{A,c}$.
\end{enumerate}
% Add Figure 4 here
%\vspace*{3cm}
Hence by Lemma~\ref{bd-solu}, we may get uniformly bounded solutions $\phi_{\varepsilon}^{A}$ and $\phi_{\varepsilon}^{B}$ of (\ref{eqn3-5-A}) and (\ref{eqn3-5-B}), respectively. Moreover, $\phi_{\varepsilon }^{A}\left( x \right)\to {{\phi }_{{{A}_{1}},0}}$ and $\phi_{\varepsilon }^{B}\left( x \right)\to {{\phi }_{{{B}_{1}},0}}$ for $x\in (-1,1)$ as $\varepsilon\to 0$. Therefore, we complete the proof of Theorem~\ref{thm1.1}.

\subsection{ Proof of Theorem~\ref{thm1.2} }\label{mul-solu-3s-2}
\ \ \ Let $N=4$, $z_2=-z_1=q_1\geq 1$ and $z_4=-z_3=q_2\geq 1$. Assume $g_{11}=g_{22}=g>0$, $g_{33}=g_{44}=\tilde{g}>0$ and ${{g}_{ij}}={{g}_{ji}}=0$ for $i=1,2$ and $j=3,4$. Then (\ref{eqn3-4}) may be represented as
\BE\label{eqn3-4-1}	
\ln {{c}_{i}}+{{z}_{i}}\phi +\sum\limits_{j=1}^{2}{{{g}_{ij}}{{c}_{j}}}=0\quad \text{ for }\;i=1,2,
\EE
and
\BE\label{eqn3-4-2}	
\ln {{c}_{i}}+{{z}_{i}}\phi +\sum\limits_{j=3}^{4}{{{g}_{ij}}{{c}_{j}}}=0\quad \text{ for }\;i=3,4\,.
\EE
Note that both (\ref{eqn3-4-1}) and (\ref{eqn3-4-2}) have the same form as (\ref{eqn-3.1}) with (\ref{eqn-3.2.1}) which can be solved explicitly. As for Theorem~\ref{pro-cn-cp} in Section~\ref{ss-solu-2s}, both (\ref{eqn3-4-1}) and (\ref{eqn3-4-2}) have two branches of solutions, respectively. We may denote these solutions as follows:
\begin{align*}
& \left({{c}_{1}},{{c}_{2}}\right)=\left({{c}_{1}}\left({{\Sigma }_{{{A}_{1}}}}\left(\phi\right)\right),{{c}_{2}}\left({{\Sigma}_{{{A}_{1}}}}\left(\phi\right)\right)\right)\,, \\
& \left({{c}_{1}},{{c}_{2}} \right)=\left({{c}_{1}}\left({{\Sigma }_{{{B}_{1}}}}\left(\phi\right) \right),{{c}_{2}}\left({{\Sigma }_{{{B}_{1}}}}\left(\phi\right)\right)\right)\,, \\
& \left({{c}_{3}},{{c}_{4}}\right)=\left({{c}_{3}}\left({{\Sigma }_{{{M}_{1}}}}\left(\phi\right)\right),{{c}_{4}}\left({{\Sigma}_{{{M}_{1}}}}\left(\phi\right)\right)\right)\,, \\
& \left({{c}_{3}},{{c}_{4}} \right)=\left({{c}_{3}}\left({{\Sigma }_{{{N}_{1}}}}\left(\phi\right) \right),{{c}_{4}}\left({{\Sigma}_{{{N}_{1}}}}\left(\phi\right)\right)\right)\,,
\end{align*}
such that
$(c_1-c_2)\circ\Sigma_{A_1}:[-\phi_{A,c},\infty)\to\mathbb{R}$, $(c_1-c_2)\circ\Sigma_{B_1}:(-\infty,\phi_{A,c}]\to\mathbb{R}$,
$(c_3-c_4)\circ\Sigma_{N_1}:[-\phi_{M,c},\infty)\to\mathbb{R}$ and $(c_3-c_4)\circ\Sigma_{M_1}:(-\infty,\phi_{M,c}]\to\mathbb{R}$, are monotone increasing functions of $\phi$, where $\phi_{A,c}, \phi_{M,c}>0$ are constants, $\Sigma_{A_1}$, $\Sigma_{B_1}$, $\Sigma_{M_1}$ and $\Sigma_{N_1}$ are functions satisfying
\begin{align*} & (c_1-c_2)\circ\Sigma_{A_1}(-\phi_{A,c})\,, (c_3-c_4) \circ\Sigma_{N_1}(-\phi_{M,c})>0\,,\\
& (c_1-c_2)\circ\Sigma_{B_1}(\phi_{A,c})\,, (c_3-c_4) \circ\Sigma_{M_1}(\phi_{M,c})<0\,, \\
& \lim_{\phi\to\infty}(c_1-c_2)\circ\Sigma_{A_1}(\phi)= \lim_{\phi\to\infty}(c_3-c_4)\circ\Sigma_{N_1}(\phi)=\infty\,, \\
& \lim_{\phi\to-\infty}(c_1-c_2)\circ\Sigma_{B_1}(\phi)= \lim_{\phi\to-\infty}(c_3-c_4)\circ\Sigma_{M_1}(\phi)=-\infty\,.
\end{align*}
Here $\circ$ denotes function composition. Moreover, Theorem~\ref{pro-cn-cp} gives $\phi_{A,c}, \phi_{M,c}\to +\infty$ and $(c_1-c_2) \circ\Sigma_{A_1}(-\phi_{A,c}), (c_1-c_2)\circ\Sigma_{B_1}(\phi_{A,c}), (c_3-c_4) \circ\Sigma_{M_1}(\phi_{M,c})$ and $(c_3-c_4) \circ\Sigma_{N_1}(-\phi_{M,c})$ tend to zero as $z, \tilde{z}\to +\infty$ and $g, \tilde{g}>0$ are fixed.

\begin{figure}[h]
 	\centering
 	\includegraphics[scale=0.37, bb=0 0 1250 180]{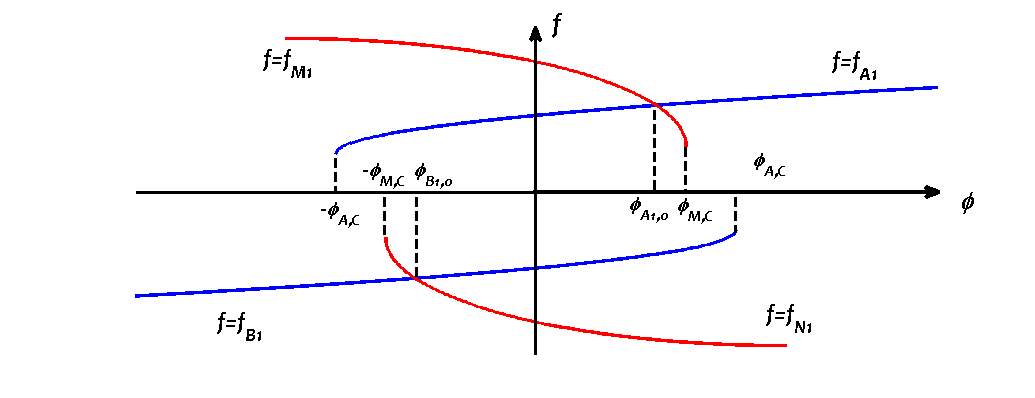}
	\vspace*{-0.7cm}  	
  	\caption{Figures of $f_{A_1}$, $f_{B_1}$, $f_{M_1}$ and $f_{N_1}$}
  	\label{fig5}
\end{figure}

Without loss of generality, we may assume $\phi_{M,c}<\phi_{A,c}$. Fix ${{\rho }_{0}}\in \mathbb{R}$ arbitrarily. Then as for (\ref{eqn-3.1}), we may solve (\ref{eqn3-4-1}) and get functions ${{f}_{{{A}_{1}}}}(\phi )=q_1\,({{c}_{1}}-{{c}_{2}})\left( {{\Sigma }_{{{A}_{1}}}}(\phi ) \right)-\rho_0$ and ${{f}_{{{B}_{1}}}}(\phi )=q_1\,({{c}_{1}}-{{c}_{2}})\left( {{\Sigma }_{{{B}_{1}}}}(\phi ) \right)-\rho_0$ which are sketched in Figure~\ref{fig5} (up to a shift by $\rho_0$), provided that $g_{12}=g_{21}=z>0$ is sufficiently large. Similarly, we may solve (\ref{eqn3-4-2}) and get functions ${{f}_{{{M}_{1}}}}(\phi )=q_2\,({{c}_{4}}-{{c}_{3}})\left( {{\Sigma }_{{{{M}}_{1}}}}(\phi ) \right)$ and ${{f}_{{{N}_{1}}}}(\phi )=q_2\,({{c}_{4}}-{{c}_{3}})\left( {{\Sigma }_{{{{N}}_{1}}}}(\phi ) \right)$ as $g_{34}=g_{43}=\tilde{z}>0$ sufficiently large (see Figure~\ref{fig5}). Because function ${{q}_{2}}({{c}_{3}}-{{c}_{4}})\circ{{\Sigma }_{{{M}_{1}}}}$ is negative and increasing to $\phi$, function ${{q}_{2}}({{c}_{4}}-{{c}_{3}})\circ {{\Sigma }_{{{M}_{1}}}}$ becomes positive and decreasing to $\phi$. On the other hand, function ${{q}_{1}}({{c}_{1}}-{{c}_{2}})\circ {{\Sigma }_{{{A}_{1}}}}-{{\rho }_{0}}$ is positive and increasing to $\phi$. This implies that as $z$ and $\tilde{z}$ sufficiently large, functions ${{q}_{1}}({{c}_{1}}-{{c}_{2}})\circ {{\Sigma }_{{{A}_{1}}}}-{{\rho }_{0}}$ and ${{q}_{2}}({{c}_{3}}-{{c}_{4}})\circ {{\Sigma }_{{{M}_{1}}}}$ may intersect at $\phi={{\phi }_{{{A}_{1}},0}}$. Similarly, functions $q_1\,(c_1-c_2)\circ\Sigma_{B_1}-{{\rho }_{0}}$ and $q_2\,(c_4-c_3) \circ\Sigma_{N_1}$ may intersect at $\phi=\phi_{B_1,0}$ as $z$ and $\tilde{z}$ sufficiently large. Generically, values $\phi_{A_1,0}$ and $\phi_{B_1,0}$ can be different by choosing $z$ and $\tilde{z}$ suitably e.g.~$z$ and $\tilde{z}$ sufficiently large.

\vspace*{2.8cm}

\noindent Let ${{f}_{A}}={{f}_{{{A}_{1}}}}-{{f}_{{{M}_{1}}}}$ and ${{f}_{B}}={{f}_{{{B}_{1}}}}-{{f}_{{{N}_{1}}}}$. Then
\[{{f}_{A}}\left( \phi  \right)={{q}_{1}}({{c}_{1}}-{{c}_{2}})\left( {{\Sigma }_{{{A}_{1}}}}(\phi ) \right)+{{q}_{2}}({{c}_{3}}-{{c}_{4}})\left( {{\Sigma }_{{{M}_{1}}}}(\phi ) \right)-{{\rho }_{0}}\,,\]
and
\[{{f}_{B}}\left( \phi  \right)={{q}_{1}}({{c}_{1}}-{{c}_{2}})\left( {{\Sigma }_{{{B}_{1}}}}(\phi ) \right)+{{q}_{2}}({{c}_{3}}-{{c}_{4}})\left( {{\Sigma }_{{{N}_{1}}}}(\phi ) \right)-{{\rho }_{0}}\,,\]
satisfy
\begin{enumerate}
\item[(1)]~~$f_A:[-\phi_{A,c},\phi_{M,c}]\to\mathbb{R}$ is smooth and strictly monotone increasing, $-\phi_{A,c}<0, f_A(-\phi_{A,c})<0, f_A(\phi_{M,c})>0$ and $f_A(\phi_{{A_1},0})=0$ for some $\phi_{{A_1},0}>-\phi_{A,c}$.
\item[(2)]~~$f_B:[-\phi_{M,c},\phi_{A,c}]\to\mathbb{R}$ is smooth and strictly monotone increasing, $\phi_{A,c}>0, f_B(\phi_{A,c})>0, f_B(-\phi_{M,c})<0$ and $f_B(\phi_{{B_1},0})=0$ for some $\phi_{{B_1},0}<\phi_{A,c}$.
\end{enumerate}
Moreover, equation (\ref{eqn-1.9}) can be expressed as $\varepsilon {{\phi }_{xx}}={{f}_{A}}\left( \phi  \right)$ and $\varepsilon {{\phi }_{xx}}={{f}_{B}}\left( \phi  \right)$ for $x\in \left( -1,1 \right)$ which have the same forms as equations~(\ref{eqn3-5-A}) and (\ref{eqn3-5-B}), respectively. Therefore by Lemma~\ref{bd-solu}, we may complete the proof of Theorem~\ref{thm1.2}.

\section{Uniformly bounded solutions}\label{sec4}

In this section, we consider the equation \BE\label{eqn-phi}
\varepsilon\phi^{\prime\prime}(x)=f(\phi(x))\quad\hbox{for}\quad
x\in (-1,1)\,, \EE with the Robin boundary condition
\BE\label{bdc-Rbn} \phi(1)+\eta_\varepsilon\phi'(1)=\phi_0(1)
\quad\hbox{and}\quad
\phi(-1)-\eta_\varepsilon\phi'(-1)=\phi_0(-1)\,,\EE where
$\phi_0(1), \phi_0(-1)$ are constants and $\eta_\varepsilon$ is a
non-negative constant. Note that the solution~$\phi_\varepsilon$ of
(\ref{eqn-phi})-(\ref{bdc-Rbn}) may depend on the
parameter~$\varepsilon$. For notational convenience, we omit $\varepsilon$ and denote $\phi$ as the solution of (\ref{eqn-phi})-(\ref{bdc-Rbn}). To get uniform boundedness of $\phi$, we assume the function $f$ satisfies one of the following conditions:
\begin{enumerate}
\item[(F1)]~~$f:[A,M]\to\mathbb{R}$ is smooth and strictly
monotone increasing, $A<0, f(A)<0, 0<M\leq\infty, f(M)>0$ and $f(\phi_A)=0$
for some $A<\phi_A<M$. \item[(F2)]~~$f:[-M,B]\to\mathbb{R}$ is
is smooth and strictly monotone increasing, $B>0, f(B)>0, 0<M\leq\infty,
f(-M)<0$ and $f(\phi_B)=0$ for some $-M<\phi_B<B$.
\end{enumerate}
Then we have

\BL\label{bd-solu} Assume the function $f$ satisfies either~(F1)
or~(F2), and the constants $A\leq\phi_0(-1)$, $\phi_0(1)\leq M$ as~(F1)
holds, and $-M\leq\phi_0(-1),\phi_0(1)\leq B$ as~(F2) holds.
Let $c=\phi_A$ if~(F1) holds, and $c=\phi_B$ if~(F2) holds. Let $\phi$
be a nonconstant solution of (\ref{eqn-phi}) with the Robin
boundary condition~(\ref{bdc-Rbn}). Then
\begin{enumerate}
\item[(i)]~~If $\phi_0(1), \phi_0(-1)>c$, then there exists
$x_1\in (-1,1)$ such that $\phi'(x_1)=0$, $\phi(x_1)>c$, and
$\phi$ is strictly monotone decreasing in $(-1,x_1)$ and
increasing in $(x_1,1)$. \item[(ii)]~~If $\phi_0(1),
\phi_0(-1)<c$, then there exists $x_2\in (-1,1)$ such that
$\phi'(x_2)=0$, $\phi(x_2)<c$, and $\phi$ is strictly monotone
increasing in $(-1,x_2)$ and decreasing in $(x_2,1)$.
\item[(iii)]~~If $\phi_0(1)\geq c\geq\phi_0(-1)$, then $\phi$ is
monotone increasing in $(-1,1)$. \item[(iv)]~~If $\phi_0(1)\leq
c\leq\phi_0(-1)$, then $\phi$ is monotone decreasing in $(-1,1)$.
\item[(v)]~~$\min\{\phi_0(-1),\phi_0(1),0\}\leq\phi(x)\leq\max\{\phi_0(-1),\phi_0(1),0\}$
for $x\in (-1,1)$. \item[(vi)]~~$\phi(x)\to c$ as $\varepsilon\to
0+$, where $c=\phi_A$ if~(F1) holds, and $c=\phi_B$ if~(F2) holds.
\item[(vii)]~~If $\ds\lim_{\varepsilon\to
0+}\frac{\varepsilon}{2\eta_\varepsilon^2}=\gamma>0$ and ${{\phi }_{0}}\left( \pm 1 \right)\neq c$, then the
solution~$\phi$ has boundary layers at $x=\pm 1$.
\end{enumerate}
\EL

\begin{proof}
Without loss of generality, we may assume the function $f$
satisfying~(F1). Replacing $\phi$ by $\phi+c$, we may assume $c=0$
and $f(0)=0$ in the whole proof for notational convenience. Since
the domain of the function~$f$ is only $[A,M]$, then we
firstly extend it smoothly to the entire real line $\mathbb{R}$ in
order to use the standard direct method to get the existence of
solution~$\phi$. Hence we may temporarily assume the function $f$
as a smooth and strictly monotone increasing function on
$\mathbb{R}$. Actually, such an assumption can be ignored because
of (\ref{phi-dom-f}).

To prove Lemma~\ref{bd-solu}, we need the following Proposition:
\BP\label{bd-solu-P}
\medskip

\noindent\begin{enumerate} \item[(a)] If $x_a\in (-1,1)$ is a
local minimum point of $\phi$, then $\phi(x_a)>0$, $\phi$ is
monotone decreasing in $(-1,x_a)$ and increasing in $(x_a,1)$.
\item[(b)] If $x_b\in (-1,1)$ is a local maximum point of $\phi$,
then $\phi(x_b)<0$, $\phi$ is monotone increasing in $(-1,x_b)$
and decreasing in $(x_b,1)$.
\end{enumerate}
\EP

The proof of Proposition~\ref{bd-solu-P}~(b) is quite similar to
that of Proposition~\ref{bd-solu-P}~(a) so we only state the proof
of Proposition~\ref{bd-solu-P}~(a) as follows: Suppose $x_a\in
(-1,1)$ is a local minimum point of $\phi$. Then $\phi'(x_a)=0$
and $\phi^{\prime\prime}(x_a)\geq 0$. If
$\phi^{\prime\prime}(x_a)=0$, then the equation
$\varepsilon\phi^{\prime\prime}=f(\phi)$ gives
$f(\phi(x_a))=\varepsilon\phi^{\prime\prime}(x_a)=0$ which implies
$\phi(x_a)=0$ and then by the uniqueness of ordinary differential
equations and $\phi(x_a)=\phi'(x_a)=0$, we have $\phi\equiv 0$
which contradicts to $\phi$ is nonconstant. Hence
$\phi^{\prime\prime}(x_a)>0$ and
$f(\phi(x_a))=\varepsilon\phi^{\prime\prime}(x_a)>0$
i.e.~$\phi(x_a)>0$. Now we prove that $\phi$ is decreasing in
$(-1,x_a)$ and increasing in $(x_a,1)$. Suppose not. Then there
exists $x_c\in (-1,1)$ and $x_c\neq x_a$ such that $x_c$ is a
local maximum point of $\phi$ i.e.~$\phi'(x_c)=0,
\phi^{\prime\prime}(x_c)\leq 0$ and $\phi(x_c)>\phi(x_1)>0$ but
$\varepsilon\phi^{\prime\prime}(x_c)=f(\phi(x_c))>0$ which
contradicts to $\phi^{\prime\prime}(x_c)\leq 0$. Therefore, we may
complete the proof of Proposition~\ref{bd-solu-P}.

For the proof~Lemma~\ref{bd-solu}~(i), we need
\medskip

\noindent {\bf Claim~1.}~~Assume $\phi_0(1), \phi_0(-1)>0$. Then
$\phi(-1), \phi(1)>0$, $\phi'(-1)<0$ and $\phi'(1)>0$. \\
We may prove Claim~1 by contradiction. Suppose one of the
following cases holds: \\
{\bf Case~I.}~~$\phi(-1)>0$ and $\phi'(-1)\geq 0$. \\
{\bf Case~II.}~~$\phi(-1)\leq 0$.

For the Case~I, we may use $\phi(-1)>0$ and the continuity of
$\phi$ to obtain that as $x\in (-1,1)$ sufficiently close to $-1$,
$\phi(x)>0$ and $\varepsilon\phi^{\prime\prime}(x)=f(\phi(x))>0$
which implies $\ds\lim_{x\to
-1+}\varepsilon\phi^{\prime\prime}(x)=f(\phi(-1))>0$. Since
$\phi'(-1)\geq 0$ and $\ds\lim_{x\to
-1+}\phi^{\prime\prime}(x)>0$, then $\phi$ is monotone increasing
in $(-1,-1+\delta_0)$, where $\delta_0>0$ is a constant. Now we
may show that $\phi$ is monotone increasing in $(-1,1)$ by
contradiction. Suppose $\phi$ has a local maximum point at $x_0\in
(-1,1)$ such that $\phi'(x_0)=0$, $\phi^{\prime\prime}(x_0)\leq 0$
and $\phi$ is monotone increasing in $(-1,x_0)$. However,
$\varepsilon\phi^{\prime\prime}(x_0)=f(\phi(x_0))\geq
f(\phi(-1))>0$ contradicts to $\phi^{\prime\prime}(x_0)\leq 0$.
Hence $\phi$ is monotone increasing in $(-1,1)$ which provides
$\phi^{\prime\prime}(x)=\frac{1}{\varepsilon}f(\phi(x))\geq\frac{1}{\varepsilon}f(\phi(-1))
$ i.e.~$\phi^{\prime\prime}(x)\geq\frac{1}{\varepsilon}f(\phi(-1))
$ for $x\in (-1,1)$. Integrating the inequality from $-1$ to $x$,
we have
$\phi'(x)-\phi'(-1)\geq\frac{1}{\varepsilon}f(\phi(-1))(x+1)$
i.e.~$\phi'(x)\geq\phi'(-1)+\frac{1}{\varepsilon}f(\phi(-1))(x+1)$
for $x\in (-1,1)$ which implies
$$\phi(1)-\phi(-1)=\int_{-1}^1\phi'(x)\,dx
\geq\int_{-1}^1\,\left[\phi'(-1)+\frac{1}{\varepsilon}f(\phi(-1))(x+1)\right]\,dx
=2\left[\phi'(-1)+\frac{1}{\varepsilon}f(\phi(-1))\right]\,,$$
i.e.~$\phi(1)\geq\phi(-1)+2\left[\phi'(-1)+\frac{1}{\varepsilon}f(\phi(-1))\right]\geq
\frac{2}{\varepsilon}f(\phi(-1))$. On the other hand, the Robin
boundary condition~(\ref{bdc-Rbn}) gives
$\phi_0(1)=\phi(1)+\eta_\varepsilon\phi'(1)\geq \phi(1)$ and
$\phi(-1)=\phi_0(-1)+\eta_\varepsilon\phi'(-1)\geq\phi_0(-1)>0$
since $\phi$ is monotone increasing in $(-1,1)$. Thus
$$\phi_0(1)\geq\phi(1)\geq\frac{2}{\varepsilon}f(\phi(-1))
\geq\frac{2}{\varepsilon}f(\phi_0(-1))\,, $$ which contradicts to
the hypothesis that $\phi_0(1), \phi_0(-1)$ are independent to
$\varepsilon$.

For the Case~II, we first use the Robin boundary
condition~(\ref{bdc-Rbn}) to get
$\eta_\varepsilon\phi'(-1)=\phi(-1)-\phi_0(-1)\leq -\phi_0(-1)<0$
which implies $\eta_\varepsilon>0$ and $\phi'(-1)<0$. Then
$\phi(x)<0$ for $x\in (-1,-1+\delta_1)$ and $\phi$ is monotone
decreasing in $(-1,-1+\delta_1)$, where $\delta_1>0$ is a
constant. Hence $\phi$ is negative and monotone decreasing in
$(-1,1)$. Otherwise, there exists $x_3\in (-1,1)$ a local minimum
point of $\phi$ such that $\phi(x_3)<0$ and
$\phi^{\prime\prime}(x_3)\geq 0$ but
$\phi^{\prime\prime}(x_3)=\frac{1}{\varepsilon}f(\phi(x_3))<0$
which contradicts to $\phi^{\prime\prime}(x_3)\geq 0$. Such a
contradiction shows that $\phi$ is negative and monotone
decreasing in $(-1,1)$. However,
$0>\phi(1)=\phi_0(1)-\eta_\varepsilon\phi'(1)\geq\phi_0(1)$
contradicts to $\phi_0(1)>0$. Notice that both Case~I and II
produce contradiction. Similarly, the condition $\phi(1)>0$ and
$\phi'(1)\leq 0$ and the other condition $\phi(1)\leq 0$ also
result in contradiction, respectively. Therefore, we may complete
the proof of Claim~I.

By Claim~I, there exists $x_1\in (-1,1)$ a local minimum point of
$\phi$, and then by Proposition~\ref{bd-solu-P}~(a), we may
complete the proof of Lemma~\ref{bd-solu}~(i). On the other hand,
we may also use the similar argument of Claim~I to prove that
there exists $x_2\in (-1,1)$ a local maximum point of $\phi$.
Hence by Proposition~\ref{bd-solu-P}~(b), we complete the proof of
Lemma~\ref{bd-solu}~(ii).

Now we prove Lemma~\ref{bd-solu}~(iii) by contradiction. Suppose
$\phi$ is not monotone increasing. By Proposition~\ref{bd-solu-P},
it is sufficient to consider two cases as follows: $\phi(-1)<0$
and $\phi(-1)>0$. If $\phi(-1)<0$, then
Proposition~\ref{bd-solu-P} implies that there exists $x_2\in
(-1,1)$ a maximum point of $\phi$ such that $\phi(x_2)<0$, $\phi$
is monotone increasing in $(-1,x_2)$ and decreasing in $(x_2,1)$
so $\phi'(1)\leq 0$. However, the boundary
condition~$\phi(1)+\eta_\varepsilon\phi'(1)=\phi_0(1)$ and
$\phi'(1)\leq 0$ give
$\phi_0(1)\leq\phi_0(1)-\eta_\varepsilon\phi'(1)=\phi(1)\leq\phi(x_2)<0$
which contradicts to $\phi_0(1)\geq c=0$. On the other hand, if
$\phi(-1)>0$, then Proposition~\ref{bd-solu-P} implies that there
exists $x_1\in (-1,1)$ a minimum point of $\phi$ such that
$\phi(x_1)>0$, $\phi$ is monotone decreasing in $(-1,x_1)$ and
increasing in $(x_1,1)$ so $\phi'(-1)\leq 0$. However, the
boundary condition~$\phi(-1)-\eta_\varepsilon\phi'(-1)=\phi_0(-1)$
and $\phi'(-1)\leq 0$ give
$\phi_0(-1)=\phi(-1)-\eta_\varepsilon\phi'(-1)\geq\phi(-1)>0$
which contradicts to $\phi_0(-1)\leq c=0$. Therefore, we complete
the proof of Lemma~\ref{bd-solu}~(iii). Similar argument of
Lemma~\ref{bd-solu}~(iii) can be applied to prove
Lemma~\ref{bd-solu}~(iv) and we omit the detail here.

Using Lemma~\ref{bd-solu}~(i)-(iv), we may prove
$\min\{\phi_0(-1),\phi_0(1),0\}\leq\phi(x)\leq\max\{\phi_0(-1),\phi_0(1),0\}$
for $x\in (-1,1)$. The proof is stated as follows: By
Lemma~\ref{bd-solu}~(i) and the boundary
condition~(\ref{bdc-Rbn}), we have
$\phi(-1)=\phi_0(-1)+\eta_\varepsilon\phi'(-1)\leq\phi_0(-1)$,
$\phi(1)=\phi_0(1)-\eta_\varepsilon\phi'(1)\leq\phi_0(1)$ and
$c=0<\phi(x_1)\leq\phi(x)\leq\max\{\phi(1),\phi(-1)\}\leq\max\{\phi_0(1),\phi_0(-1)\}$
for $x\in (-1,1)$. Similarly, Lemma~\ref{bd-solu}~(ii) and the
boundary condition~(\ref{bdc-Rbn}) imply
$\phi(-1)=\phi_0(-1)+\eta_\varepsilon\phi'(-1)\geq\phi_0(-1)$,
$\phi(1)=\phi_0(1)-\eta_\varepsilon\phi'(1)\geq\phi_0(1)$ and
$c=0>\phi(x_2)\geq\phi(x)\geq\min\{\phi(1),\phi(-1)\}\geq\min\{\phi_0(1),\phi_0(-1)\}$
for $x\in (-1,1)$. On the other hand, we may apply
Lemma~\ref{bd-solu} (iii) and the boundary
condition~(\ref{bdc-Rbn}) to get
$\phi(-1)=\phi_0(-1)+\eta_\varepsilon\phi'(-1)\geq\phi_0(-1)$,
$\phi(1)=\phi_0(1)-\eta_\varepsilon\phi'(1)\leq\phi_0(1)$ and
$\phi_0(-1)\leq\phi(-1)\leq\phi(x)\leq\phi(1)\leq\phi_0(1)$ for
$x\in (-1,1)$. Similarly, Lemma~\ref{bd-solu}~(iv) and the
boundary condition~(\ref{bdc-Rbn}) give
$\phi(-1)=\phi_0(-1)+\eta_\varepsilon\phi'(-1)\leq\phi_0(-1)$,
$\phi(1)=\phi_0(1)-\eta_\varepsilon\phi'(1)\geq\phi_0(1)$ and
$\phi_0(-1)\geq\phi(-1)\geq\phi(x)\geq\phi(1)\geq\phi_0(1)$ for
$x\in (-1,1)$. Hence we complete the proof of
Lemma~\ref{bd-solu}~(v) i.e. \BE\label{max-norm-phi}
\min\{\phi_0(-1),\phi_0(1),0\}\leq\phi(x)\leq\max\{\phi_0(-1),\phi_0(1),0\}\quad\hbox{
for}\quad x\in (-1,1)\,.\EE Let
$A_0=\min\{\phi_0(-1),\phi_0(1),0\}$ and
$A_1=\max\{\phi_0(-1),\phi_0(1),0\}$. Then (\ref{max-norm-phi})
implies \BE\label{infty-norm-phi} \|\phi\|_{L^\infty}\leq
A_2=\max\{-A_0, A_1\}\,. \EE Since $A\leq\phi_0(-1), \phi_0(1)\leq M$
and $A<0$, then (\ref{max-norm-phi}) gives
\BE\label{phi-dom-f}\phi(x)\in [A_0,A_1]\subset
[A,M]\quad\hbox{ for}\quad x\in (-1,1)\,, \EE i.e.~each value
of $\phi(x)$ must be contained in the original domain of the
function $f$. Thus we may neglect the extension of the
function~$f$ and regard $\phi$ as a well-defined solution of
equation~(\ref{eqn-phi}) with boundary
condition~(\ref{bdc-Rbn}).

Now we claim that $\phi(x)\to 0$ as $\varepsilon\to 0+$ for $x\in
(-1,1)$. To prove this, we remark that
$$\frac{\varepsilon}{2}(\phi^2)^{\prime\prime}(x)=\varepsilon(\phi\phi^{\prime\prime}+(\phi')^2)(x)
\geq\varepsilon\phi\phi^{\prime\prime}(x)=\phi(x)\,f(\phi(x))=\phi(x)\,\int_0^{\phi(x)}\,f'(s)\,ds\geq
\alpha_0\phi^2(x)\,, $$ for $x\in (-1,1)$, where
$\alpha_0=\min_{z\in [A_0,A_1]}f'(z)>0$ is a constant coming from
the strictly monotone increasing of $f$. Note that if $\phi(x)<0$,
then
$$\phi(x)\,\int_0^{\phi(x)}\,f'(s)\,ds=(-\phi(x))\,\int_{\phi(x)}^0\,f'(s)\,ds
\geq(-\phi(x))\,\int_{\phi(x)}^0\,\alpha_0\,ds=\alpha_0\phi^2(x)\,.
$$ Since
$\frac{\varepsilon}{2}(\phi^2)^{\prime\prime}(x)\geq\alpha_0\phi^2(x)$
for $x\in (-1,1)$, then by (\ref{infty-norm-phi}) and the standard
comparison theorem, we have $\phi^2(x)\leq
A_2^2\left(e^{-(1+x)\sqrt{2\alpha_0/\varepsilon}}+e^{-(1-x)\sqrt{2\alpha_0/\varepsilon}}\right)$
for $x\in (-1,1)$. Therefore, $\phi(x)\to 0$ as $\varepsilon\to
0+$ for $x\in (-1,1)$, and we may complete the proof of
Lemma~\ref{bd-solu}~(vi).

For the proof of Lemma~~\ref{bd-solu}~(vii), we firstly multiply
the equation~(\ref{eqn-phi}) by $\phi'$. Then we have $
\frac{\varepsilon}{2}\left[(\phi')^2\right]'=\varepsilon\phi'\phi^{\prime\prime}=f(\phi)\phi'=\frac{d}{dx}F(\phi)$,
which implies \BE\label{phi'1}
\frac{\varepsilon}{2}(\phi'(x))^2=F(\phi(x))+C_\varepsilon\,, \EE
for $x\in (-1,1)$, where $F(\phi)=\ds\int_0^\phi f(s)\,ds$ and
$C_\varepsilon$ is a constant depending on $\varepsilon$. Now we
want to claim that $\ds\lim_{\varepsilon\to 0+}\,C_\varepsilon=0$.
By the mean value theorem, there exists $x_\varepsilon\in
(-\frac{1}{2},\frac{1}{2})$ such that
$\phi(\frac{1}{2})-\phi(-\frac{1}{2})=\phi'(x_\varepsilon)$. Since
$\phi(x)\to 0$ as $\varepsilon\to 0+$ for $x\in (-1,1)$, then
$\phi(\frac{1}{2}), \phi(-\frac{1}{2})\to 0$ and
$\phi'(x_\varepsilon)=\phi(\frac{1}{2})-\phi(-\frac{1}{2})\to 0$
as $\varepsilon\to 0+$. Hence by (\ref{phi'1}), we obtain
$C_\varepsilon=\frac{\varepsilon}{2}(\phi'(x_\varepsilon))^2-F(\phi(x_\varepsilon))\to
0$ as $\varepsilon\to 0+$ i.e.~$\ds\lim_{\varepsilon\to
0+}\,C_\varepsilon=0$. On the other hand, we may put the Robin
boundary condition~(\ref{bdc-Rbn}) into (\ref{phi'1}) and get
$\frac{\varepsilon}{2\eta_\varepsilon^2}\left(\phi_0(1)-\phi(1)\right)^2=F(\phi(1))+C_\varepsilon
$ and
$\frac{\varepsilon}{2\eta_\varepsilon^2}\left(\phi_0(-1)-\phi(-1)\right)^2=F(\phi(-1))+C_\varepsilon
$. By (\ref{phi-dom-f}) and the continuity of $\phi$, we may
assume $\phi(\pm 1)\to \phi_*(\pm 1)$ as $\varepsilon\to 0+$ (up
to a subsequence). Generically, the values $\phi_*(1)$ and
$\phi_*(-1)$ may not be uniquely determined but here we want to
claim the uniqueness of $\phi_*(\pm 1)$ as follows: Suppose
$\ds\lim_{\varepsilon\to
0+}\frac{\varepsilon}{2\eta_\varepsilon^2}=\gamma>0$. Then
$\phi_*(1)$ and $\phi_*(-1)$ satisfy
$\sqrt{\gamma}|\phi_0(1)-\phi_*(1)|=\sqrt{F(\phi_*(1))}$ and
$\sqrt{\gamma}|\phi_0(-1)-\phi_*(-1)|=\sqrt{F(\phi_*(-1))}$.
Notice that the function $F$ is positive and monotone increasing
in $(0,M]$ and decreasing in $[A,0)$ because ${F}'\left( \phi  \right)=f\left( \phi  \right)>0$ on $\left( 0,M \right)$ and $<0$ on $\left( A,0 \right)$. Here we have used the fact that ${{\phi }_{A}}=c=0$. By
Lemma~\ref{bd-solu}~(i)-(iv), we have $\sqrt{\gamma}|\phi_0(\pm
1)-\phi_*(\pm 1)|=\sqrt{\gamma}(\phi_0(\pm 1)-\phi_*(\pm 1))$ if
$\phi_0(\pm 1)>0$; $\sqrt{\gamma}|\phi_0(\pm 1)-\phi_*(\pm
1)|=\sqrt{\gamma}(\phi_*(\pm 1)-\phi_0(\pm 1))$ if $\phi_0(\pm
1)<0$; $\sqrt{\gamma}|\phi_0(\pm 1)-\phi_*(\pm
1)|=\pm\sqrt{\gamma}(\phi_0(\pm 1)-\phi_*(\pm 1))$ if $\phi_0(
-1)\leq 0\leq\phi_0(1)$; and $\sqrt{\gamma}|\phi_0(\pm
1)-\phi_*(\pm 1)|=\pm\sqrt{\gamma}(\phi_*(\pm 1)-\phi_0(\pm 1))$
if $\phi_0(-1)\geq 0\geq\phi_0(1)$. Hence $\phi_*(\pm 1)$ can be
uniquely determined by the equations
$\sqrt{\gamma}|\phi_0(1)-s|=\sqrt{F(s)}$ and
$\sqrt{\gamma}|\phi_0(-1)-s|=\sqrt{F(s)}$, respectively. The
uniqueness of $\phi_*(\pm 1)$ implies that the asymptotic limits
of boundary values $\phi(\pm 1)$ are $\ds\lim_{\varepsilon\to
0+}\phi(\pm 1)= \phi_*(\pm 1)$. One may also remark that
$\phi_*(\pm 1)\neq 0=c$ (Otherwise, if $\phi_*(\pm 1)=0$,then $\sqrt{\gamma}|\phi_0(1)-\phi_*(1)|=\sqrt{F(\phi_*(1))}$,
$\sqrt{\gamma}|\phi_0(-1)-\phi_*(-1)|=\sqrt{F(\phi_*(-1))}$ and $F(0)=0$
imply that ${{\phi }_{0}}\left( \pm 1 \right)=0=c$ which contradicts to the assumption ${{\phi }_{0}}\left( \pm 1 \right)\neq c$.) Consequently, the solution~$\phi$ has
boundary layers at $x=\pm 1$ if $\gamma>0$. Therefore, we have
showed Lemma~\ref{bd-solu}~(vii) and completed the proof of
Lemma~\ref{bd-solu}.
\end{proof}

\BR\label{uniq-phi} The equation~(\ref{eqn-phi}) with the boundary
condition~(\ref{bdc-Rbn}) has a unique solution. \ER

The uniqueness comes from the strictly monotone increasing of the
function~$f$. The proof is sketched as follows: Suppose $\phi_1$
and $\phi_2$ are solutions of (\ref{eqn-phi}) and~(\ref{bdc-Rbn}).
We may subtract the equation of $\phi_1$ by that of $\phi_2$, and
multiply the resulting equation by $u=\phi_1-\phi_2$ and integrate
it over $(-1,1)$. Then using integration by part, we have
$u'(1)u(1)-u'(-1)u(-1)-\int_{-1}^1\,(u'(x))^2\,dx=\int_{-1}^1\,c(x)u^2\,dx
$, where
$c(x)=\frac{f(\phi_1(x))-f(\phi_2(x))}{\phi_1(x)-\phi_2(x)}$ is
positive since the function $f$ is strictly monotone increasing.
On the other hand, the Robin boundary condition~(\ref{bdc-Rbn})
gives $u(-1)=\eta_\varepsilon u'(-1)$, $u(1)=-\eta_\varepsilon
u'(1)$ and
$u'(1)u(1)-u'(-1)u(-1)=-\eta_\varepsilon\left[(u'(-1))^2+(u'(1))^2
\right]$. Hence
$$0\leq\int_{-1}^1\,c(x)u^2\,dx=-\eta_\varepsilon\left[(u'(-1))^2+(u'(1))^2
\right]-\int_{-1}^1\,(u'(x))^2\,dx\leq 0$$ which implies $u\equiv
0$ i.e. $\phi_1\equiv\phi_2$ and the uniqueness proof of $\phi$ is
complete.

\BR\label{lin-stab-phi} The solution~$\phi$ of the
equation~(\ref{eqn-phi}) with the boundary
condition~(\ref{bdc-Rbn}) has linear stability. \ER

To get the linear stability of the solution~$\phi$ of the
equation~(\ref{eqn-phi}) with the boundary
condition~(\ref{bdc-Rbn}), we study the eigenvalue problem
$Lv=\lambda v$ of the corresponding linearized operator
$Lv=-\varepsilon v^{\prime\prime}+f'(\phi)v$ with the boundary
condition $v(\pm 1)\pm\eta_\varepsilon v(\pm 1)=0$. Using
integration by part, it is obvious that
\begin{align}\notag &\lambda\int_{-1}^1 v^2\,dx=\int_{-1}^1
v\,Lv\,dx=\int_{-1}^1 \varepsilon v^{\prime\prime}
v\,dx+\int_{-1}^1 f'(\phi)\,v^2\,dx
\\ &=\eta_\varepsilon\left[(v'(-1))^2+(v'(1))^2\right]
+\int_{-1}^1 \left[\varepsilon(v')^2+f'(\phi)v^2\right]\,dx\geq
\mu_0\int_{-1}^1 v^2\,dx\,, \end{align} and hence
$\lambda\geq\mu_0>0$, where $\mu_0=\min_{s\in
[\min\phi,\max\phi]}f'(s)$ is a positive constant arising from the
strictly monotone increasing of the function~$f$.

In Lemma~\ref{bd-solu}, the existence of zero point $\phi_A$ (or
$\phi_B$) of $f$ is essential. If the function~$f$ has not any
zero point like $\phi_A$ (or $\phi_B$) i.e. the value of~$f$ is
away from zero, then the equation~(\ref{eqn-phi}) may not have
uniformly bounded solutions $\{\phi\}_{\varepsilon>0}$. Such a
result is stated as follows:
\BL\label{ub-solu} Assume $f$ is a
function satisfying one of the following conditions:
\begin{enumerate} \item[(a)]~~$f:[A,\infty)\to\mathbb{R}$ is
monotone increasing, $A<0$ and $f(A)>0$.
\item[(b)]~~$f:(-\infty,B]\to\mathbb{R}$ is monotone increasing,
$B>0$ and $f(B)<0$.
\end{enumerate}
For each $\varepsilon>0$, let $\phi$ be a solution of the
equation~(\ref{eqn-phi}).
Then~$\ds\sup_{\varepsilon>0}\,\|\phi\|_{L^\infty}=\infty$. \EL
\medskip
\noindent

\begin{proof}
Without loss of generality, we may assume the function $f$
satisfies the condition~(a). Now we prove Lemma~\ref{ub-solu} by
contradiction. Suppose~$\{\phi\}_{\varepsilon>0}$ is uniformly
bounded
i.e.~$\ds\sup_{\varepsilon>0}\,\|\phi\|_{L^\infty}<\infty$. We
divide three cases to complete the proof as follows:
\medskip

\noindent {\bf Case~I.}~~The solution $\phi=\phi(x)$ is monotone
decreasing to $x$ i.e.~$\phi'(x)\leq 0$ for $x\in (-1,1)$.

Using the equation $\varepsilon\phi^{\prime\prime}=f(\phi) $ and
the condition~(a), we have
\begin{align} -\phi'(x)\geq
\phi'(1)-\phi'(x)&=\int_x^1\,\phi^{\prime\prime}(\tau)\,d\tau \notag \\
&=\varepsilon^{-1}\int_x^1\,f(\phi(\tau))\,d\tau \notag \\
&\geq\varepsilon^{-1}\int_x^1\,f(A)\,d\tau=\varepsilon^{-1}\,f(A)(1-x)\,,\quad\forall
x\in (-1,1)\,, \notag
\end{align} and hence
$$-2\|\phi\|_{L^\infty}\leq\phi(1)-\phi(-1)=\ds\int_{-1}^1\phi'(x)\,dx\leq -\varepsilon^{-1}f(A)\int_{-1}^1\,(1-x)\,dx
=-2\varepsilon^{-1}f(A)\,, $$
i.e.~$\|\phi\|_{L^\infty}\geq\varepsilon^{-1}f(A)\to\infty$ as
$\varepsilon\to 0+$ which contradicts to the
hypothesis~$\ds\sup_{\varepsilon>0}\,\|\phi\|_{L^\infty}<\infty$.
\medskip

\noindent {\bf Case~II.}~~The solution $\phi=\phi(x)$ is monotone
increasing to $x$ i.e.~$\phi'(x)\geq 0$ for $x\in (-1,1)$.

As for the argument of Case~I, we obtain
\begin{align} \phi'(x)\geq
\phi'(x)-\phi'(-1)&=\int_{-1}^x\,\phi^{\prime\prime}(\tau)\,d\tau \notag \\
&=\varepsilon^{-1}\int_{-1}^x\,f(\phi(\tau))\,d\tau \notag \\
&\geq\varepsilon^{-1}\int_{-1}^x\,f(A)\,d\tau=\varepsilon^{-1}\,f(A)(1+x)\,,\quad\forall
x\in (-1,1)\,, \notag
\end{align} and hence
$$2\|\phi\|_{L^\infty}\geq\phi(1)-\phi(-1)=\ds\int_{-1}^1\phi'(x)\,dx\geq \varepsilon^{-1}f(A)\int_{-1}^1\,(1+x)\,dx
=2\varepsilon^{-1}f(A)\,, $$
i.e.~$\|\phi\|_{L^\infty}\geq\varepsilon^{-1}f(A)\to\infty$ as
$\varepsilon\to 0+$ which contradicts to the
hypothesis~$\ds\sup_{\varepsilon>0}\,\|\phi\|_{L^\infty}<\infty$.
\medskip

\noindent {\bf Case~III.}~~The solution $\phi=\phi(x)$ has a local
minimum point at $x_0\in (-1,1)$ such that~$\phi'(x_0)=0$ and
$\phi^{\prime\prime}(x_0)>0$.

Note that since $\varepsilon\phi^{\prime\prime}=f(\phi)\geq
f(A)>0$, it is impossible to have any local maximum point in
$(-1,1)$. By the equation $\varepsilon\phi^{\prime\prime}=f(\phi)
$ and the condition~(a), we have
\begin{align} -\phi'(x)=
\phi'(x_0)-\phi'(x)&=\int_x^{x_0}\,\phi^{\prime\prime}(\tau)\,d\tau \notag \\
&=\varepsilon^{-1}\int_x^{x_0}\,f(\phi(\tau))\,d\tau \notag \\
&\geq\varepsilon^{-1}\int_x^{x_0}\,f(A)\,d\tau=\varepsilon^{-1}\,f(A)(x_0-x)\,,\quad\forall
x\in (-1,x_0)\,, \notag
\end{align} and hence
$$-2\|\phi\|_{L^\infty}\leq\phi(x_0)-\phi(-1)=\ds\int_{-1}^{x_0}\phi'(x)\,dx\leq -\varepsilon^{-1}f(A)\int_{-1}^{x_0}\,(x_0-x)\,dx
=-\frac{1}{2}\varepsilon^{-1}f(A)(x_0+1)^2\,, $$ i.e.
\BE\label{x0-1} |x_0+1|\leq
2\varepsilon^{1/2}\,\sqrt{\|\phi\|_{L^\infty}/f(A)}\,. \EE

On the other hand, \begin{align} \phi'(x)=
\phi'(x)-\phi'(x_0)&=\int_{x_0}^x\,\phi^{\prime\prime}(\tau)\,d\tau \notag \\
&=\varepsilon^{-1}\int_{x_0}^x\,f(\phi(\tau))\,d\tau \notag \\
&\geq\varepsilon^{-1}\int_{x_0}^x\,f(A)\,d\tau=\varepsilon^{-1}\,f(A)(x-x_0)\,,\quad\forall
x\in (x_0,1)\,, \notag
\end{align} and hence
$$2\|\phi\|_{L^\infty}\geq\phi(1)-\phi(x_0)=\ds\int_{x_0}^1\phi'(x)\,dx\geq \varepsilon^{-1}f(A)\int_{x_0}^1\,(x-x_0)\,dx
=\frac{1}{2}\varepsilon^{-1}f(A)(x_0-1)^2\,, $$ i.e.
\BE\label{x0-2} |x_0-1|\leq
2\varepsilon^{1/2}\,\sqrt{\|\phi\|_{L^\infty}/f(A)}\,. \EE
Therefore, as $\varepsilon>0$ sufficiently small, (\ref{x0-1}) and
(\ref{x0-2}) provide a contradiction and we may complete the proof
of Lemma~\ref{ub-solu}.
\end{proof}

\section{Excess currents due to steric effects}\label{ex-cu}
Here we want to use solutions $\phi_{\varepsilon}^{A}$ and $\phi_{\varepsilon}^{B}$ of (\ref{eqn3-4})-(\ref{bdc-Rbn-0}) (see Theorem~\ref{thm1.1} and~\ref{thm1.2}) to calculate excess currents (due to steric effects) represented by formula~(\ref{I-ex1}). By (\ref{eqn3-1}), $$\sum\limits_{j=1}^{N}{{{g}_{ij}}{{c}_{j}}=-{{k}_{B}}T\ln {{c}_{i}}-{{z}_{i}}e\phi }\,,$$ and then formula~(\ref{I-ex1}) becomes
\BE\label{I-ex1.1}
{{I}^{ex}}=\sum\limits_{i=1}^{N}{{{z}_{i}}e{{D}_{i}}\left( \nabla {{c}_{i}}+{{z}_{i}}{{c}_{i}}\nabla \tilde{\phi } \right)}\,,
\EE
where $\tilde{\phi }=\frac{e}{{{k}_{B}}T}\phi $.
\subsection{Under the same hypotheses of Theorem~\ref{thm1.1}}\label{ex-cu1}
\ \ \  Here we set $N=3$, $z_2=-z_1=q\geq 1, z_3>0$, $\rho_0>0$, and assume that $g_{11}=g_{22}=g>0$ is fixed, $g_{12}=g_{21}=z>0$ is sufficiently large, and $g_{i3}=g_{3i}=0$ for $i=1,2,3$. By (\ref{eqn-3.2}), we have ${{c}_{3}}={{e}^{-{{z}_{3}}\tilde{\phi }}}$ which implies $\nabla {{c}_{3}}+{{z}_{3}}{{c}_{3}}\nabla \tilde{\phi }=0$. Hence (\ref{I-ex1.1}) becomes
\begin{eqnarray*} {{I}^{ex}}&=&\sum\limits_{i=1}^{2}{{{z}_{i}}e{{D}_{i}}\left( \frac{d{{c}_{i}}}{dx}+{{z}_{i}}{{c}_{i}}\frac{d\tilde{\phi }}{dx} \right)} \\
&=&-q\,e{{D}_{1}}\left( \frac{d{{c}_{1}}}{dx}-q\,{{c}_{1}}\frac{d\tilde{\phi }}{dx} \right)+q\,e{{D}_{2}}\left( \frac{d{{c}_{2}}}{dx}+q\,{{c}_{2}}\frac{d\tilde{\phi }}{dx} \right) \\
&=& q\,e\left[ \frac{d}{dx}\left( -{{D}_{1}}{{c}_{1}}+{{D}_{2}}{{c}_{2}} \right)+q\,\left( {{D}_{1}}{{c}_{1}}+{{D}_{2}}{{c}_{2}} \right)\frac{d\tilde{\phi }}{dx} \right]
\end{eqnarray*}	
i.e. \BE\label{I-ex1.2}
{{I}^{ex}}=q\,e\left[ \frac{d}{dx}\left( -{{D}_{1}}{{c}_{1}}+{{D}_{2}}{{c}_{2}} \right)+q\,\left( {{D}_{1}}{{c}_{1}}+{{D}_{2}}{{c}_{2}} \right)\frac{d\tilde{\phi }}{dx} \right]
\EE
Using  ${{c}_{1}}=\frac{{{c}_{1}}-{{c}_{2}}}{2}+\frac{{{c}_{1}}+{{c}_{2}}}{2}$ and ${{c}_{2}}=\frac{{{c}_{1}}+{{c}_{2}}}{2}-\frac{{{c}_{1}}-{{c}_{2}}}{2}$, formula (\ref{I-ex1.2}) can be expressed as
\BE\label{I-ex1.3}
\begin{array}{rll}\frac{1}{q\,e}\,{{I}^{ex}}=&\frac{d}{dx}\left[ \frac{{{D}_{2}}-{{D}_{1}}}{2}\left( {{c}_{1}}+{{c}_{2}} \right)-\frac{{{D}_{1}}+{{D}_{2}}}{2}\left( {{c}_{1}}-{{c}_{2}} \right) \right] \\
& \\
&+q\,\left[ \frac{{{D}_{1}}+{{D}_{2}}}{2}\left( {{c}_{1}}+{{c}_{2}} \right)-\frac{{{D}_{2}}-{{D}_{1}}}{2}\left( {{c}_{1}}-{{c}_{2}} \right) \right]\frac{d\tilde{\phi }}{dx}
\end{array}
\EE
Note that ${{c}_{1}}+{{c}_{2}}=\Sigma $ and ${{c}_{1}}-{{c}_{2}}=\left\{ \begin{array}{rrrll}
   {} & \sqrt{{{\Sigma }^{2}}-4{{e}^{-\left( g+z \right)\Sigma }}} & \text{ on }\;A,  \\
   {} & -\sqrt{{{\Sigma }^{2}}-4{{e}^{-\left( g+z \right)\Sigma }}} & \text{on}\;B.  \\
\end{array} \right. $ (see (\ref{cn-cp}) in Section~\ref{ss-solu-2s}). As for (\ref{multi-sol1})-(\ref{eqn3-5-B}), we may set $\left( \Sigma ,\tilde{\phi } \right)=\left({{\Sigma }_{{{A}_{1}}}}\left( {\tilde{\phi }}\right),\phi_{\varepsilon }^{A}\left( x \right) \right)$ and $\left( \Sigma ,\tilde{\phi } \right)=\left({{\Sigma }_{{{B}_{1}}}}\left( {\tilde{\phi }} \right),\phi_{\varepsilon }^{B}\left( x \right) \right)$ , respectively. Then along $c_1+c_2=\Sigma ={{\Sigma }_{{{A}_{1}}}}$ and $c_1-c_2=\sqrt{{\Sigma^2_{A_1}}-4{{e}^{-\left( g+z \right)\Sigma_{A_1}}}}$, we may use (\ref{d-rho-a+}), (\ref{d-cn-cp-a+}) and Chain Rule to get
\[\begin{array}{*{35}{l}}
   \frac{d}{dx}\left( {{c}_{1}}+{{c}_{2}} \right) & =\frac{d}{dx}{{\Sigma }_{{{A}_{1}}}}\left( \phi _{\varepsilon }^{A}\left( x \right) \right)  \\
   {} & =\frac{d{{\Sigma }_{{{A}_{1}}}}}{d\phi }\left( \phi _{\varepsilon }^{A}\left( x \right) \right)\frac{d\phi _{\varepsilon }^{A}}{dx}\left( x \right)  \\
   {} & =\frac{\sqrt{\Sigma_{{{A}_{1}}}^{2}\left( \phi _{\varepsilon }^{A}\left( x \right) \right)-4{{e}^{-(g+z){{\Sigma }_{{{A}_{1}}}}\left( \phi _{\varepsilon }^{A}\left( x \right) \right)}}}}{1+g{{\Sigma }_{{{A}_{1}}}}\left( \phi _{\varepsilon }^{A}\left( x \right) \right)+\left( {{g}^{2}}-{{z}^{2}} \right){{e}^{-(g+z){{\Sigma }_{{{A}_{1}}}}\left( \phi _{\varepsilon }^{A}\left( x \right) \right)}}}\frac{d\phi _{\varepsilon }^{A}}{dx}\left( x \right)\,,  \\
\end{array}\]
and	
\[\begin{array}{*{35}{l}}
   \frac{d}{dx}\left( {{c}_{1}}-{{c}_{2}} \right) & =\frac{d}{dx}\left( {{c}_{1}}-{{c}_{2}} \right)\left( {{\Sigma }_{{{A}_{1}}}}\left( \phi _{\varepsilon }^{A}\left( x \right) \right) \right)  \\
   {} & =\frac{d}{d\phi }\left( {{c}_{1}}-{{c}_{2}} \right)\left( {{\Sigma }_{{{A}_{1}}}}\left( \phi _{\varepsilon }^{A}\left( x \right) \right) \right)\frac{d\phi _{\varepsilon }^{A}}{dx}\left( x \right)  \\
   {} & =\frac{{{\Sigma }_{{{A}_{1}}}}\left( \phi _{\varepsilon }^{A}\left( x \right) \right)+2(g+z){{e}^{-(g+z){{\Sigma }_{{{A}_{1}}}}\left( \phi _{\varepsilon }^{A}\left( x \right) \right)}}}{1+g{{\Sigma }_{{{A}_{1}}}}\left( \phi _{\varepsilon }^{A}\left( x \right) \right)+({{g}^{2}}-{{z}^{2}}){{e}^{-(g+z){{\Sigma }_{{{A}_{1}}}}\left( \phi _{\varepsilon }^{A}\left( x \right) \right)}}}\frac{d\phi_{\varepsilon }^{A}}{dx}\left( x \right)\,.  \\
\end{array}\]
For simplicity, we may set ${{\hat{\Sigma }}_{{{A}_{1}}}}={{\Sigma}_{{{A}_{1}}}}\left(\phi_{\varepsilon }^{A}\left( x \right)\right)$
and denote $\frac{d}{dx}\left( {{c}_{1}}\pm {{c}_{2}} \right)$ as follows:
\[\frac{d}{dx}\left( {{c}_{1}}+{{c}_{2}} \right)=\frac{\sqrt{\hat{\Sigma }_{{{A}_{1}}}^{2}-4{{e}^{-(g+z){{{\hat{\Sigma }}}_{{{A}_{1}}}}}}}}{1+g{{{\hat{\Sigma }}}_{{{A}_{1}}}}+\left( {{g}^{2}}-{{z}^{2}} \right){{e}^{-(g+z){{{\hat{\Sigma }}}_{{{A}_{1}}}}}}}\frac{d\phi _{\varepsilon }^{A}}{dx}\left( x \right)\,,\]
and
\[\frac{d}{dx}\left( {{c}_{1}}-{{c}_{2}} \right)=\frac{{{{\hat{\Sigma }}}_{{{A}_{1}}}}+2(g+z){{e}^{-(g+z){{{\hat{\Sigma }}}_{{{A}_{1}}}}}}}{1+g{{{\hat{\Sigma }}}_{{{A}_{1}}}}+({{g}^{2}}-{{z}^{2}}){{e}^{-(g+z){{{\hat{\Sigma }}}_{{{A}_{1}}}}}}}\frac{d\phi _{\varepsilon }^{A}}{dx}\left( x \right)\,.\]
Consequently, by setting $I_A^{ex}=I^{ex}$ along $A_1$, (\ref{I-ex1.3}) becomes
\BE\label{I-ex1.4}
\begin{array}{*{35}{l}}
\frac{1}{q\,e}{{I}_A^{ex}}=&\frac{\frac{{{D}_{2}}-{{D}_{1}}}{2}\sqrt{\hat{\Sigma }_{{{A}_{1}}}^{2}-4{{e}^{-\left( g+z \right){{{\hat{\Sigma }}}_{{{A}_{1}}}}}}}-\frac{{{D}_{1}}+{{D}_{2}}}{2}\left[ {{{\hat{\Sigma }}}_{{{A}_{1}}}}+2\left( g+z \right){{e}^{-\left( g+z \right){{{\hat{\Sigma }}}_{{{A}_{1}}}}}} \right]}{1+g{{{\hat{\Sigma }}}_{{{A}_{1}}}}+\left( {{g}^{2}}-{{z}^{2}} \right){{e}^{-\left( g+z \right){{{\hat{\Sigma }}}_{{{A}_{1}}}}}}}\frac{d\phi _{\varepsilon }^{A}}{dx}\left( x \right)  \\
& \\
& +q\,\left[\frac{{{D}_{1}}+{{D}_{2}}}{2}{{{\hat{\Sigma }}}_{{{A}_{1}}}} - \frac{{{D}_{2}}-{{D}_{1}}}{2}\sqrt{\hat{\Sigma }_{{{A}_{1}}}^{2}-4{{e}^{-\left( g+z\right){{{\hat{\Sigma }}}_{{{A}_{1}}}}}}}\right]\frac{d\phi _{\varepsilon }^{A}}{dx}\left( x \right)  \\
\end{array}
\EE
Similarly, along ${{c}_{1}}+{{c}_{2}}=\Sigma ={{\Sigma }_{{{B}_{1}}}}\left(\tilde{\phi}\right)$, $c_1-c_2=-\sqrt{{\Sigma^2_{B_1}}-4{{e}^{-\left( g+z \right)\Sigma_{B_1}}}}$ and $\tilde{\phi }=\phi_{\varepsilon }^{B}\left( x \right)$, we may set $I^{ex}=I_B^{ex}$, and use (\ref{d-rho-b+}) and (\ref{d-cn-cp-b+}) to get
\BE\label{I-ex1.5}
\begin{array}{*{35}{l}}
\frac{1}{q\,e}{{I}_B^{ex}}=&\frac{\frac{{{D}_{1}}-{{D}_{2}}}{2}\sqrt{\hat{\Sigma }_{{{B}_{1}}}^{2}-4{{e}^{-\left( g+z \right){{{\hat{\Sigma }}}_{{{B}_{1}}}}}}}-\frac{{{D}_{1}}+{{D}_{2}}}{2}\left[ {{{\hat{\Sigma }}}_{{{B}_{1}}}}+2\left( g+z \right){{e}^{-\left( g+z \right){{{\hat{\Sigma }}}_{{{B}_{1}}}}}} \right]}{1+g{{{\hat{\Sigma }}}_{{{B}_{1}}}}+\left( {{g}^{2}}-{{z}^{2}} \right){{e}^{-\left( g+z \right){{{\hat{\Sigma }}}_{{{B}_{1}}}}}}}\frac{d\phi _{\varepsilon }^{B}}{dx}\left( x \right)  \\
& \\
& +q\,\left[\frac{{{D}_{1}}+{{D}_{2}}}{2}{{{\hat{\Sigma }}}_{{{B}_{1}}}} - \frac{{{D}_{1}}-{{D}_{2}}}{2}\sqrt{\hat{\Sigma }_{{{B}_{1}}}^{2}-4{{e}^{-\left( g+z\right){{{\hat{\Sigma }}}_{{{B}_{1}}}}}}}\right]\frac{d\phi _{\varepsilon }^{B}}{dx}\left( x \right)  \\
\end{array}
\EE
where ${{\hat{\Sigma }}_{{{B}_{1}}}}={{\Sigma }_{{{B}_{1}}}}\left( \phi _{\varepsilon }^{B}\left( x \right) \right)$.

Equation~(\ref{I-ex1.4}) and (\ref{I-ex1.5}) can be denoted as
\BE\label{I-ex1.4.1}
I_{A}^{ex}=q\,e\,{{i}_{A}}\left( {{{\hat{\Sigma }}}_{{{A}_{1}}}} \right)\frac{d\phi^{A}_{\varepsilon }}{dx}\left( x \right)\,,
\EE
and
\BE\label{I-ex1.5.1}
I_{B}^{ex}=q\,e\,{{i}_{B}}\left( {{{\hat{\Sigma }}}_{{{B}_{1}}}} \right)\frac{d\phi^{B}_{\varepsilon }}{dx}\left( x \right)\,,
\EE
where
\BE\label{iA1}
\begin{array}{*{35}{l}}
{{i}_{A}}\left( \Sigma  \right)=&\frac{\frac{{{D}_{2}}-{{D}_{1}}}{2}\sqrt{{{\Sigma }^{2}}-4{{e}^{-\left( g+z \right)\Sigma }}}-\frac{{{D}_{1}}+{{D}_{2}}}{2}\left[ \Sigma +2\left( g+z \right){{e}^{-\left( g+z \right)\Sigma }} \right]}{1+g\Sigma +\left( {{g}^{2}}-{{z}^{2}} \right){{e}^{-\left( g+z \right)\Sigma }}} \,, \\
& \\
& +q\,\left[\frac{{{D}_{1}}+{{D}_{2}}}{2}\Sigma-\frac{{{D}_{2}}-{{D}_{1}}}{2}\sqrt{{{\Sigma }^{2}}-4{{e}^{-\left( g+z \right)\Sigma }}} \right]  \\
\end{array}
\EE
and
\BE\label{iB1}
\begin{array}{*{35}{l}}
{{i}_{B}}\left( \Sigma  \right)= &\frac{\frac{{{D}_{1}}-{{D}_{2}}}{2}\sqrt{{{\Sigma }^{2}}-4{{e}^{-\left( g+z \right)\Sigma }}}-\frac{{{D}_{1}}+{{D}_{2}}}{2}\left[ \Sigma +2\left( g+z \right){{e}^{-\left( g+z \right)\Sigma }} \right]}{1+g\Sigma +\left( {{g}^{2}}-{{z}^{2}} \right){{e}^{-\left( g+z \right)\Sigma }}}\,.  \\
& \\
& +q\,\left[\frac{{{D}_{1}}+{{D}_{2}}}{2}\Sigma-\frac{{{D}_{1}}-{{D}_{2}}}{2}\sqrt{{{\Sigma }^{2}}-4{{e}^{-\left( g+z \right)\Sigma }}} \right]  \\
\end{array}
\EE
Without loss of generality, $\phi_{\varepsilon}^{A}$ can be assumed as a monotone increasing function. Such an assumption can be fulfilled by setting $\phi_0(-1)<\phi_0(1)$ and using Lemma~\ref{bd-solu}~(iii). Integrating $I_A^{ex}$ from $x_1$ to $x_2$, we have
\BE\label{iA-2}
\int_{{{x}_{1}}}^{{{x}_{2}}}{I_{A}^{ex}dx}=e\,\int_{{{x}_{1}}}^{{{x}_{2}}}{{{i}_{A}}\left( {{\Sigma }_{{{A}_{1}}}}\left( \phi _{\varepsilon }^{A}\left( x \right) \right) \right)}\frac{d\phi _{\varepsilon }^{A}}{dx}dx=e\,\int_{{{\phi }^A_{1}}}^{{{\phi}^A_{2}}}{{{i}_{A}}\left( {{\Sigma }_{{{A}_{1}}}}\left( \phi  \right) \right)}d\phi\,,
\EE
for $-1<x_1<x_2<1$, where $\phi_1\leq\phi_2$ and ${{\phi}^A_{j}}=\phi _{\varepsilon }^{A}\left( {{x}_{j}} \right), j=1,2$. Setting $\Sigma ={{\Sigma }_{{{A}_{1}}}}$ and using change of variables, Inverse Function Theorem and (\ref{d-rho-a+}), we have
\[d\phi =\frac{d\phi }{d\Sigma }d\Sigma =\frac{1}{\frac{d\Sigma }{d\phi }}d\Sigma =\frac{1+g\Sigma +\left( {{g}^{2}}-{{z}^{2}} \right){{e}^{-(g+z)\Sigma }}}{\sqrt{{{\Sigma }^{2}}-4{{e}^{-(g+z)\Sigma }}}}d\Sigma\,.\]
Then by (\ref{iA1}), we get
\begin{align} \nonumber
  & {{i}_{A}}\left( \Sigma  \right)d\phi ={{i}_{A}}\left( \Sigma  \right)\frac{1+g\Sigma +\left( {{g}^{2}}-{{z}^{2}} \right){{e}^{-(g+z)\Sigma }}}{\sqrt{{{\Sigma }^{2}}-4{{e}^{-(g+z)\Sigma }}}}d\Sigma  \\ \nonumber
 & =\frac{{{D}_{2}}-{{D}_{1}}}{2}\left\{(1-q)-q\left[ g\Sigma +\left( {{g}^{2}}-{{z}^{2}} \right){{e}^{-(g+z)\Sigma }} \right]\right\}d\Sigma \\ \nonumber
 & \hspace{0.3cm} -\frac{{{D}_{1}}+{{D}_{2}}}{2\sqrt{{{\Sigma }^{2}}-4{{e}^{-(g+z)\Sigma }}}}\left[(1-q)\Sigma+ 2\left( g+z \right){{e}^{-(g+z)\Sigma }}-q\,g{{\Sigma }^{2}}-q\,\left( {{g}^{2}}-{{z}^{2}} \right)\Sigma {{e}^{-(g+z)\Sigma }} \right]d\Sigma  \\ \nonumber
\end{align}
Hence (\ref{iA-2}) becomes
\BE\label{iA-3}
\begin{array}{rll}
&\ds\int_{{{x}_{1}}}^{{{x}_{2}}}{I_{A}^{ex}dx}= q\,e\,\ds\int_{{{\phi}^A_{1}}}^{{{\phi}^A_{2}}}\,{{{i}_{A}}\left( {{\Sigma }_{{{A}_{1}}}}\left( \phi  \right) \right)}d\phi \\
& \\
&=q\,e\,\ds\int_{{{\Sigma}^A_{1}}}^{{{\Sigma}^A_{2}}}\,\frac{{{D}_{2}}-{{D}_{1}}}{2}\left\{(1-q)-q\left[ g\Sigma +\left( {{g}^{2}}-{{z}^{2}} \right){{e}^{-(g+z)\Sigma }} \right]\right\}d\Sigma \\
& \\
& \hspace{0.3cm} -q\,e\,\ds\int_{{{\Sigma}^A_{1}}}^{{{\Sigma}^A_{2}}}\,\frac{{{D}_{1}}+{{D}_{2}}}{2\sqrt{{{\Sigma }^{2}}-4{{e}^{-(g+z)\Sigma }}}}\left\{(1-q)\Sigma-q\,g{{\Sigma }^{2}}+(g+z)\left[2-q\,\left( {{g}}-{{z}} \right)\Sigma\right] {{e}^{-(g+z)\Sigma }} \right\}d\Sigma\,,
\end{array}
\EE
where ${{\Sigma}^A_{j}}={{\Sigma }_{{{A}_{1}}}}\left( {{\phi}^A_{j}} \right)$ and $\phi_{j}^{A}= \left(\phi_{\varepsilon }^{A}\left( {{x}_{j}} \right) \right)$  for $j=1,2$. Similarly, we may set $\Sigma ={{\Sigma }_{{{B}_{1}}}}$ and use change of variables, Inverse Function Theorem and (\ref{d-rho-b+}) to get
\[d\phi =\frac{d\phi }{d\Sigma }d\Sigma =\frac{1}{\frac{d\Sigma }{d\phi }}d\Sigma =-\frac{1+g\Sigma +\left( {{g}^{2}}-{{z}^{2}} \right){{e}^{-(g+z)\Sigma }}}{\sqrt{{{\Sigma }^{2}}-4{{e}^{-(g+z)\Sigma }}}}d\Sigma\,.\]
Then as for (\ref{iA-3}), we may use (\ref{iB1}) to derive
\BE\label{iB-3}
\begin{array}{rll}
&\ds\int_{{{x}_{1}}}^{{{x}_{2}}}{I_{B}^{ex}dx}=q\,e\,\ds\int_{{{\phi}^B_{1}}}^{{{\phi}^B_{2}}}{}{{i}_{B}}\left( {{\Sigma }_{{{B}_{1}}}}\left( \phi\right) \right)d\phi \\ & \\ &=q\,e\,\ds\int_{{{\Sigma}^B_{1}}}^{{{\Sigma}^B_{2}}}{}
\frac{{{D}_{2}}-{{D}_{1}}}{2}\left\{(1-q)-q\left[ g\Sigma +\left( {{g}^{2}}-{{z}^{2}} \right){{e}^{-(g+z)\Sigma }} \right]\right\}d\Sigma \\
& \\ & \hspace{0.3cm} +q\,e\,\ds\int_{{{\Sigma}^B_{1}}}^{{{\Sigma }^B_{2}}}{}\frac{{{D}_{1}}+{{D}_{2}}}{2\sqrt{{{\Sigma }^{2}}-4{{e}^{-(g+z)\Sigma }}}}\left\{(1-q)\Sigma-q\,g{{\Sigma }^{2}}+(g+z)\left[2-q\,\left( {{g}}-{{z}} \right)\Sigma\right] {{e}^{-(g+z)\Sigma }} \right\}d\Sigma\,,
\end{array}
\EE
where ${{\Sigma}^B_{j}}={{\Sigma }_{{{B}_{1}}}}\left( {{\phi}^B_{j}} \right)$ and $\phi_{j}^{B}=\phi_{\varepsilon }^{B}\left( {{x}_{j}} \right)$ for $j=1,2$. Therefore, we complete the proof of (\ref{iA-3-int}) and (\ref{iB-3-int}).

\subsection{Under the same hypotheses of Theorem~\ref{thm1.2}}\label{ex-cu2}
\ \ \ Here we set $N=4$, $z_2=-z_1=q_1\geq 1$, $z_4=-z_3=q_2\geq 1$, $\rho_0\neq 0$, and assume that $g_{11}=g_{22}=g>0$, $g_{33}=g_{44}=\tilde{g}>0$ are fixed, $g_{12}=g_{21}=z>0$, $g_{34}=g_{43}=\tilde{z}>0$ are sufficiently large, and ${{g}_{ij}}={{g}_{ji}}=0$ for $i=1,2$ and $j=3,4$. As for Section~\ref{mul-solu-3s-2}, these hypotheses imply that (\ref{eqn3-4}) can be decomposed into two independent equations (\ref{eqn3-4-1}) and (\ref{eqn3-4-2}) which have the same form as (\ref{eqn-3.1}) with (\ref{eqn-3.2.1}). Solving equations (\ref{eqn3-4-1}) and (\ref{eqn3-4-2}), we may get $(c_1,c_2)$ (with branches $A_1, B_1$) and $(c_3,c_4)$ (with branches $M_1,N_1$) as functions of $\phi$, respectively. By (\ref{I-ex1.1}), the excess currents of $(c_1,c_2)$ and $(c_3,c_4)$ can be represented as $\sum\limits_{i=1}^{2}{{{z}_{i}}e{{D}_{i}}\left( \frac{d{{c}_{i}}}{dx}+{{z}_{i}}{{c}_{i}}\frac{d\tilde{\phi }}{dx} \right)}$ and $\sum\limits_{i=3}^{4}{{{z}_{i}}e{{D}_{i}}\left( \frac{d{{c}_{i}}}{dx}+{{z}_{i}}{{c}_{i}}\frac{d\tilde{\phi }}{dx} \right)}$ which can be calculated by the same method as Section~\ref{ex-cu1}. We may denote the total excess current as $I^{ex}_{A,M}=I^{ex}_A+I^{ex}_M$, where $I^{ex}_A$ and $I^{ex}_M$ are the excess currents along branches $A_1$ for $(c_1,c_2)$ and $M_1$ for $(c_3,c_4)$, respectively. Then as for (\ref{I-ex1.4}), we have
\BE\label{I-ex1.4A}
\begin{array}{*{35}{l}}
\frac{1}{q_1\,e}{{I}_A^{ex}}=&\frac{\frac{{{D}_{2}}-{{D}_{1}}}{2}\sqrt{\hat{\Sigma }_{{{A}_{1}}}^{2}-4{{e}^{-\left( g+z \right){{{\hat{\Sigma }}}_{{{A}_{1}}}}}}}-\frac{{{D}_{1}}+{{D}_{2}}}{2}\left[ {{{\hat{\Sigma }}}_{{{A}_{1}}}}+2\left( g+z \right){{e}^{-\left( g+z \right){{{\hat{\Sigma }}}_{{{A}_{1}}}}}} \right]}{1+g{{{\hat{\Sigma }}}_{{{A}_{1}}}}+\left( {{g}^{2}}-{{z}^{2}} \right){{e}^{-\left( g+z \right){{{\hat{\Sigma }}}_{{{A}_{1}}}}}}}\frac{d\phi _{\varepsilon }^{A}}{dx}\left( x \right)  \\
& \\
& +q_1\,\left[\frac{{{D}_{1}}+{{D}_{2}}}{2}{{{\hat{\Sigma }}}_{{{A}_{1}}}} - \frac{{{D}_{2}}-{{D}_{1}}}{2}\sqrt{\hat{\Sigma }_{{{A}_{1}}}^{2}-4{{e}^{-\left( g+z\right){{{\hat{\Sigma }}}_{{{A}_{1}}}}}}}\right]\frac{d\phi _{\varepsilon }^{A}}{dx}\left( x \right)  \\
\end{array}
\EE
and
\BE\label{I-ex1.4M}
\begin{array}{*{35}{l}}
\frac{1}{q_2\,e}{{I}_M^{ex}}=&\frac{\frac{{{D}_{4}}-{{D}_{3}}}{2}\sqrt{\hat{\Sigma }_{{{M}_{1}}}^{2}-4{{e}^{-\left( \tilde{g}+\tilde{z} \right){{{\hat{\Sigma }}}_{{{M}_{1}}}}}}}-\frac{{{D}_{3}}+{{D}_{4}}}{2}\left[ {{{\hat{\Sigma }}}_{{{M}_{1}}}}+2\left( \tilde{g}+\tilde{z} \right){{e}^{-\left( \tilde{g}+\tilde{z} \right){{{\hat{\Sigma }}}_{{{M}_{1}}}}}} \right]}{1+\tilde{g}{{{\hat{\Sigma }}}_{{{M}_{1}}}}+\left( {{\tilde{g}}^{2}}-{{\tilde{z}}^{2}} \right){{e}^{-\left( \tilde{g}+\tilde{z} \right){{{\hat{\Sigma }}}_{{{M}_{1}}}}}}}\frac{d\phi _{\varepsilon }^{A}}{dx}\left( x \right)  \\
& \\
& +q_2\,\left[\frac{{{D}_{3}}+{{D}_{4}}}{2}{{{\hat{\Sigma }}}_{{{M}_{1}}}} - \frac{{{D}_{4}}-{{D}_{3}}}{2}\sqrt{\hat{\Sigma }_{{{M}_{1}}}^{2}-4{{e}^{-\left( \tilde{g}+\tilde{z}\right){{{\hat{\Sigma }}}_{{{M}_{1}}}}}}}\right]\frac{d\phi _{\varepsilon }^{A}}{dx}\left( x \right)  \\
\end{array}
\EE
Similarly, another total excess current can be denoted as $I^{ex}_{B,N}=I^{ex}_B+I^{ex}_N$, where $I^{ex}_B$ and $I^{ex}_N$ are the excess currents along branches $B_1$ for $(c_1,c_2)$ and $N_1$ for $(c_3,c_4)$, respectively. As for (\ref{I-ex1.5}), we have
\BE\label{I-ex1.5B}
\begin{array}{*{35}{l}}
\frac{1}{q_1\,e}{{I}_B^{ex}}=&\frac{\frac{{{D}_{1}}-{{D}_{2}}}{2}\sqrt{\hat{\Sigma }_{{{B}_{1}}}^{2}-4{{e}^{-\left( g+z \right){{{\hat{\Sigma }}}_{{{B}_{1}}}}}}}-\frac{{{D}_{1}}+{{D}_{2}}}{2}\left[ {{{\hat{\Sigma }}}_{{{B}_{1}}}}+2\left( g+z \right){{e}^{-\left( g+z \right){{{\hat{\Sigma }}}_{{{B}_{1}}}}}} \right]}{1+g{{{\hat{\Sigma }}}_{{{B}_{1}}}}+\left( {{g}^{2}}-{{z}^{2}} \right){{e}^{-\left( g+z \right){{{\hat{\Sigma }}}_{{{B}_{1}}}}}}}\frac{d\phi _{\varepsilon }^{B}}{dx}\left( x \right)  \\
& \\
& +q_1\,\left[\frac{{{D}_{1}}+{{D}_{2}}}{2}{{{\hat{\Sigma }}}_{{{B}_{1}}}} - \frac{{{D}_{1}}-{{D}_{2}}}{2}\sqrt{\hat{\Sigma }_{{{B}_{1}}}^{2}-4{{e}^{-\left( g+z\right){{{\hat{\Sigma }}}_{{{B}_{1}}}}}}}\right]\frac{d\phi _{\varepsilon }^{B}}{dx}\left( x \right)  \\
\end{array}
\EE
and
\BE\label{I-ex1.4N}
\begin{array}{*{35}{l}}
\frac{1}{q_2\,e}{{I}_N^{ex}}=&\frac{\frac{{{D}_{4}}-{{D}_{3}}}{2}\sqrt{\hat{\Sigma }_{{{N}_{1}}}^{2}-4{{e}^{-\left( \tilde{g}+\tilde{z} \right){{{\hat{\Sigma }}}_{{{N}_{1}}}}}}}-\frac{{{D}_{3}}+{{D}_{4}}}{2}\left[ {{{\hat{\Sigma }}}_{{{N}_{1}}}}+2\left( \tilde{g}+\tilde{z} \right){{e}^{-\left( \tilde{g}+\tilde{z} \right){{{\hat{\Sigma }}}_{{{N}_{1}}}}}} \right]}{1+\tilde{g}{{{\hat{\Sigma }}}_{{{N}_{1}}}}+\left( {{\tilde{g}}^{2}}-{{\tilde{z}}^{2}} \right){{e}^{-\left( \tilde{g}+\tilde{z} \right){{{\hat{\Sigma }}}_{{{N}_{1}}}}}}}\frac{d\phi _{\varepsilon }^{B}}{dx}\left( x \right)  \\
& \\
& +q_2\,\left[\frac{{{D}_{3}}+{{D}_{4}}}{2}{{{\hat{\Sigma }}}_{{{N}_{1}}}} - \frac{{{D}_{4}}-{{D}_{3}}}{2}\sqrt{\hat{\Sigma }_{{{N}_{1}}}^{2}-4{{e}^{-\left( \tilde{g}+\tilde{z}\right){{{\hat{\Sigma }}}_{{{N}_{1}}}}}}}\right]\frac{d\phi _{\varepsilon }^{B}}{dx}\left( x \right)  \\
\end{array}
\EE
Hence as for (\ref{iA-3}) and (\ref{iB-3}), we may use (\ref{I-ex1.4A})-(\ref{I-ex1.4N}) to get
\BE\label{iA-3A}
\begin{array}{rll}
&\ds\int_{{{x}_{1}}}^{{{x}_{2}}}{I_{A}^{ex}dx} \\
&=q_1\,e\,\ds\int_{{{\Sigma}^A_{1}}}^{{{\Sigma}^A_{2}}}\,\frac{{{D}_{2}}-{{D}_{1}}}{2}\left\{(1-q_1)-q_1\left[ g\Sigma +\left( {{g}^{2}}-{{z}^{2}} \right){{e}^{-(g+z)\Sigma }} \right]\right\}d\Sigma \\
& \\
& \hspace{0.3cm} -q_1\,e\,\ds\int_{{{\Sigma}^A_{1}}}^{{{\Sigma}^A_{2}}}\,\frac{{{D}_{1}}+{{D}_{2}}}{2\sqrt{{{\Sigma }^{2}}-4{{e}^{-(g+z)\Sigma }}}}\left\{(1-q_1)\Sigma-q_1\,g{{\Sigma }^{2}}+(g+z)\left[2-q_1\,\left( {{g}}-{{z}} \right)\Sigma\right] {{e}^{-(g+z)\Sigma }} \right\}d\Sigma\,,
\end{array}
\EE
\BE\label{iA-3M}
\begin{array}{rll}
&\ds\int_{{{x}_{1}}}^{{{x}_{2}}}{I_{M}^{ex}dx} \\
&=q_2\,e\,\ds\int_{{{\Sigma}^M_{1}}}^{{{\Sigma}^M_{2}}}\,\frac{{{D}_{4}}-{{D}_{3}}}{2}\left\{(1-q_2)-q_2\left[ \tilde{g}\Sigma +\left( {{\tilde{g}}^{2}}-{{\tilde{z}}^{2}} \right){{e}^{-(\tilde{g}+\tilde{z})\Sigma }} \right]\right\}d\Sigma \\
& \\
& \hspace{0.3cm} -q_2\,e\,\ds\int_{{{\Sigma}^M_{1}}}^{{{\Sigma}^M_{2}}}\,\frac{{{D}_{3}}+{{D}_{4}}}{2\sqrt{{{\Sigma }^{2}}-4{{e}^{-(\tilde{g}+\tilde{z})\Sigma }}}}\left\{(1-q_2)\Sigma-q_2\,\tilde{g}{{\Sigma }^{2}}+(\tilde{g}+\tilde{z})\left[2-q_2\,\left( {{\tilde{g}}}-{{\tilde{z}}} \right)\Sigma\right] {{e}^{-(\tilde{g}+\tilde{z})\Sigma }} \right\}d\Sigma\,,
\end{array}
\EE
\BE\label{iB-3B}
\begin{array}{rll}
&\ds\int_{{{x}_{1}}}^{{{x}_{2}}}{I_{B}^{ex}dx} \\ &=q_1\,e\,\ds\int_{{{\Sigma}^B_{1}}}^{{{\Sigma}^B_{2}}}{}
\frac{{{D}_{2}}-{{D}_{1}}}{2}\left\{(1-q_1)-q_1\left[ g\Sigma +\left( {{g}^{2}}-{{z}^{2}} \right){{e}^{-(g+z)\Sigma }} \right]\right\}d\Sigma \\
& \\ & \hspace{0.3cm} +q_1\,e\,\ds\int_{{{\Sigma}^B_{1}}}^{{{\Sigma }^B_{2}}}{}\frac{{{D}_{1}}+{{D}_{2}}}{2\sqrt{{{\Sigma }^{2}}-4{{e}^{-(g+z)\Sigma }}}}\left\{(1-q_1)\Sigma-q_1\,g{{\Sigma }^{2}}+(g+z)\left[2-q_1\,\left( {{g}}-{{z}} \right)\Sigma\right] {{e}^{-(g+z)\Sigma }} \right\}d\Sigma\,,
\end{array}
\EE
\BE\label{iB-3N}
\begin{array}{rll}
&\ds\int_{{{x}_{1}}}^{{{x}_{2}}}{I_{N}^{ex}dx} \\
&=q_2\,e\,\ds\int_{{{\Sigma}^N_{1}}}^{{{\Sigma}^N_{2}}}\,\frac{{{D}_{4}}-{{D}_{3}}}{2}\left\{(1-q_2)-q_2\left[ \tilde{g}\Sigma +\left( {{\tilde{g}}^{2}}-{{\tilde{z}}^{2}} \right){{e}^{-(\tilde{g}+\tilde{z})\Sigma }} \right]\right\}d\Sigma \\
& \\
& \hspace{0.3cm} +q_2\,e\,\ds\int_{{{\Sigma}^N_{1}}}^{{{\Sigma}^N_{2}}}\,\frac{{{D}_{3}}+{{D}_{4}}}{2\sqrt{{{\Sigma }^{2}}-4{{e}^{-(\tilde{g}+\tilde{z})\Sigma }}}}\left\{(1-q_2)\Sigma-q_2\,\tilde{g}{{\Sigma }^{2}}+(\tilde{g}+\tilde{z})\left[2-q_2\,\left( {{\tilde{g}}}-{{\tilde{z}}} \right)\Sigma\right] {{e}^{-(\tilde{g}+\tilde{z})\Sigma }} \right\}d\Sigma\,,
\end{array}
\EE
where ${{\Sigma}^A_{j}}={{\Sigma }_{{{A}_{1}}}}\left(\phi_{\varepsilon }^{A}\left( {{x}_{j}} \right) \right)$, ${{\Sigma}^M_{j}}={{\Sigma }_{{{M}_{1}}}}\left(\phi_{\varepsilon }^{A}\left( {{x}_{j}} \right) \right)$, ${{\Sigma}^B_{j}}={{\Sigma }_{{{B}_{1}}}}\left(\phi_{\varepsilon }^{B}\left( {{x}_{j}} \right) \right)$, and ${{\Sigma}^N_{j}}={{\Sigma }_{{{N}_{1}}}}\left(\phi_{\varepsilon }^{B}\left( {{x}_{j}} \right) \right)$  for $j=1,2$.
Combining (\ref{iA-3A})-(\ref{iB-3N}) and using $I_{A, M}^{ex}=I_{A}^{ex}+I_{M}^{ex}, I_{B, N}^{ex}=I_{B}^{ex}+I_{N}^{ex}$, we may complete the proof of (\ref{iA-3AM}) and (\ref{iB-3BN}).

%------------------------------------------------------------------------------------------------------------

\section{Acknowledgements}
Tai-Chia Lin is partially supported by the National Science Council of Taiwan under grant numbers NSC-102-2115-M-002-015 and NSC-100-2115-M-002-007. Bob Eisenberg is supported in part by the Bard Endowed Chair of Rush University.

%% References
%-----------------------------------------------------------------------------------------------------------
%\bibliographystyle{plain}
%\bibliography{multi-PNP}

%-----------------------------------------------------------------------------------------------------------

\end{document}